\documentclass[10pt]{article}
\usepackage{psfrag,amsmath}
\usepackage{amssymb}
\usepackage{mathrsfs}
\usepackage{amsfonts}
\usepackage{amsmath}
\usepackage{mathrsfs,color}
\usepackage{epstopdf}
\usepackage{graphicx}
\usepackage{mathrsfs,amscd,amssymb,amsthm,amsmath,bm,url}
\usepackage{tikz}
\usepackage{listings}
\usepackage{titlesec}

\makeatletter

\renewcommand{\@seccntformat}[1]{{\csname the#1\endcsname}{\normalsize .}\hspace{.5em}}
\makeatother
\renewcommand{\thesection}{\normalsize \arabic{section}}

\def \[{\begin{equation}}
\def \]{\end{equation}}

\def \det{{\rm det}}
\newtheorem{thm}{Theorem}[section]

\newtheorem{claim}{Claim}

\newtheorem{lem}[thm]{Lemma}
\newtheorem{cor}[thm]{Corollary}

\newenvironment{wst}
{\setlength{\leftmargini}{1.5\parindent}
 \begin{itemize}
 \setlength{\itemsep}{-1.1mm}}
{\end{itemize}}
\begin{document}

\setlength{\baselineskip}{15pt}
\begin{center}{\Large \bf On the $A_\alpha$-index of graphs with given order and dissociation number\footnote{S.L.\ acknowledges the financial support from the National Natural Science Foundation of China (Grant Nos. 12171190,\,11671164).}}
\vspace{4mm}

{\large Zihan Zhou,\ \ \ Shuchao Li\footnote{{\it Email addresses}: 
lscmath@ccnu.edu.cn (S. Li)}}\vspace{2mm}

{\small Faculty of Mathematics and Statistics, Key Lab NAA--MOE, and Hubei Key Lab--Math. Sci.,\\ Central China Normal University, Wuhan 430079, China}
\end{center}

\noindent {\bf Abstract}:\
Given a graph $G,$ a subset of vertices is called a maximum dissociation set of $G$ if it induces a subgraph with vertex degree at most 1, and the subset has maximum cardinality. The cardinality of a maximum dissociation set is called the dissociation number of $G$. The adjacency matrix and the degree diagonal matrix of $G$ are denoted by $A(G)$ and $D(G),$ respectively. In 2017, Nikiforov proposed the $A_\alpha$-matrix: $A_\alpha(G)=\alpha D(G)+(1-\alpha)A(G),$ where $\alpha\in[0,1].$ The largest eigenvalue of this novel matrix is called the $A_\alpha$-index of $G.$  In this paper, we firstly determine the connected graph (resp. bipartite graph, tree) having the largest $A_\alpha$-index over all connected graphs (resp. bipartite graphs, trees) with fixed order and dissociation number. Secondly, we describe the structure of all the $n$-vertex graphs having the minimum $A_\alpha$-index with dissociation number $\tau$, where $\tau\geqslant\lceil\frac{2}{3}n\rceil.$ Finally, we identify all the connected $n$-vertex graphs with dissociation number $\tau\in\{2,\lceil\frac{2}{3}n\rceil,n-1,n-2\}$ having the minimum $A_\alpha$-index.

\vspace{2mm} \noindent{\it Keywords:}
$A_\alpha$-index; Dissociation number; Bipartite graph; Tree
\vspace{2mm}


\section{\normalsize Introduction}\setcounter{equation}{0}
In this section, we introduce some definitions that will help the reader to build up the necessary background for the main results.
\subsection{\normalsize Background and definitions}

Given a graph $G,$ its \textit{adjacency matrix} $A(G)$ is an $n\times n$ 0-1 square matrix whose $(u,v)$-entry is 1 if and only if $u\sim v$ in $G.$ Let $D(G)=diag(d_1,\ldots,d_n)$ be the diagonal matrix of vertex degrees in a graph $G.$ The matrix $Q(G)=D(G)+A(G)$ is called the \textit{signless Laplacian matrix} of $G;$ see \cite{D-R}. Nikiforov \cite{V.N} introduced the \textit{$A_\alpha$-matrix} of a graph $G,$ which is a convex combination of $D(G)$ and $A(G),$ i.e.,
$$
A_\alpha(G)=\alpha D(G)+(1-\alpha)A(G),\,\,\ 0\leqslant \alpha\leqslant 1.
$$
Note that $A_\alpha(G)$ is real symmetric. Hence its eigenvalues are real, and so we may display the eigenvalues of $A_\alpha(G)$ as $\lambda_1(G)\geqslant \lambda_2(G)\geqslant\cdots\geqslant \lambda_n(G).$ For short, the $A_\alpha$ spectral radius of $G$ (i.e. the largest eigenvalue of $A_\alpha(G)$), denoted by $\lambda_\alpha(G)=\lambda_1(G)$, is called the \textit{$A_\alpha$-index} of $G.$ Notice that
\begin{align}\label{eq:1.1}
A(G)=A_0(G), \ \ Q(G)=2A_{1/2}(G) \ \ \text{and} \ \ D(G)=A_1(G).
\end{align}

Recently, more and more people studied the $A_\alpha$-spectra of graphs. Nikiforov et al. \cite{C.Y} gave some bounds on the $A_\alpha$-index of a graph, and they determined the unique tree with maximum (resp. minimum) $A_\alpha$-index among $n$-vertex trees. Nikiforov and Rojo \cite{N-O} determined the graph of order $n$ and diameter at least $d$ having the largest $A_\alpha$-index. Xue et al. \cite{HZ} determined the graphs with the maximum (resp. minimum) $A_\alpha$-index among all connected graphs with given diameter (resp. clique number). Li and Zhou \cite{L-Z} determined the unique $n$-vertex block graph with prescribed independence number having the maximum $A_\alpha$-index. In the same paper, they also characterized the unique $n$-vertex graph with given number of cut edges having the largest $A_\alpha$-index, and it is surprising to see that in both cases, the extremal graphs always coincide. For more advances on the $A_\alpha$-spectra, we refer the reader to \cite{CLM,HLX,L-S,LTYZ,L-W,LY-2,LY-1,LHX,WWT,XWT}.

For a subset $S\subseteq V_G,$ we call $S$ a \textit{dissociation set} if the induced subgraph $G[S]$ does not contain $P_3$ as a subgraph, which is equivalent to say that  the degree of each vertex in $G[S]$ is at most one. A \textit{maximum dissociation set} of $G$ is a dissociation set with the maximum cardinality. The \textit{dissociation number} of $G$, written as $\tau(G)$, is the cardinality of a maximum dissociation set of $G.$ The problem of determining $\tau(G)$ is a classical problem, which may date back to the work of Yannakakis \cite{Y81} in 1981. In the same paper, he also proved that this problem is NP-complete for bipartite graphs. Boliac, Cameron and Lozin \cite{BCL} strengthened the result of Yannakakis by reducing the problem, in polynomial time, from general bipartite graphs to some particular classes such as bipartite graphs with maximum degree 3 or $C_4$-free bipartite graphs. In the same paper, they also proved that finding the dissociation number is polynomially solvable for bipartite graphs containing no induced subgraph isomorphic to a tree obtained from $P_4$ by attaching a path of length three to one quasi-pendant vertex of $P_4.$ For more advances about this problem one may consult \cite{ABKL,BPPR,BPPR1,CH,OD,PY}.

Quite recently, characterizing the graphs having maximum number of maximum dissociation sets becomes an attractive problem. Tu, Zhang and Shi \cite{TZS} determined all the trees having the maximum number of maximum dissociation sets among trees with given order. Tu, Zhang and Du \cite{TZD} characterized all the trees having the maximum number of maximum dissociation sets among trees with given dissociation number. Li and Sun \cite{LS} identified all the trees (resp. forests) having the largest and the second largest number of maximum dissociation sets among trees (resp. forests) with given order and dissociation number. Tu, Li and Du \cite{TLD} presented the upper bounds on the number of maximal (resp. maximum) dissociation sets in a general graph of order $n$ and in a triangle-free graph of order $n,$ and they also characterized the corresponding extremal graphs.

Naturally, one may be interested in relating the spectra of a graph and its dissociation number. Very recently, Huang, Li and Zhou \cite{HLZ} characterized all the graphs among all connected graphs (resp. bipartite graphs, trees) with given order and dissociation number having the maximum $A_0$-spectral radius. Das and Mohanty \cite{DM2023} identified the unique block graph with given order and dissociation number having the largest $A_0$-spectral radius.

In this paper we consider the relation between the $A_\alpha$-spectra of a graph and its dissociation number, which extends the main results of \cite{HLZ}. We firstly characterize the $n$-vertex connected graph (resp. bipartite graph, tree) with dissociation number $\tau$ having the largest $A_\alpha$-index. Secondly, we describe the structure of all the $n$-vertex graphs with dissociation number $\tau$ having the minimum $A_\alpha$-index, where $\tau\geqslant\left\lceil\frac{2}{3}n\right\rceil.$ Finally, we identify all the connected $n$-vertex graphs with dissociation number $\tau\in\{2,\left\lceil\frac{2}{3}n\right\rceil,n-1,n-2\}$ having the minimum $A_\alpha$-index.

\subsection{\normalsize Basic notations and main results}
In this subsection, we give some basic notations and then describe our main results. Throughout this paper, we consider only simple and finite graphs. For spectral graph theoretic notation and terminology not defined here, we refer to Godsil and Royle \cite{G-R} and West \cite{West}.

Let $G=(V_G,E_G)$ be a graph with vertex set $V_G$ and edge set $E_G.$ We call $n=|V_G|$ and $m=|E_G|$ the \textit{order} and the \textit{size} of $G,$ respectively.  We say that two vertices $u$ and $v$ are \textit{adjacent} (or \textit{neighbors}) if they are joined by an edge and we write $u\sim v.$  As usual, let $P_n,C_n, K_n$ and $K_{t,n-t}$ denote the path, cycle, complete graph and complete bipartite graph on $n$ vertices, respectively. Denote by $K_{1,n-1}$ the star graph, which is also denoted by $S_n$ as usual.

For a vertex $v\in V_G,$ let $N_G(v)$ be the set of all neighbors of $v$ in $G.$ Then $d_G(v)=|N_G(v)|$ is the \textit{degree} of $v$ in $G.$ For simplicity, when there is no danger of confusion, we may omit the subscripts $G$ for our notation. The vertex with degree $n-1$ in a star $S_n$ is called the \textit{centre} of $S_n.$ A \textit{pendant vertex} (or a \textit{leaf}) is a vertex of $G$ whose degree is one. We denote by $\mathcal{P}_G$ the set of all pendant vertices of $G.$ We call a vertex \textit{quasi-pendant vertex} of $G$ if it is adjacent to a pendant vertex of $G.$ We denote by $\mathcal{Q}_G$ the set of all quasi-pendant vertices of $G$ and let $\mathcal{Q}_G'$ consist of all quasi-pendant vertices of degree 2 of $G.$ A \textit{pendant path} of length $r$ is a path with vertex set $V=\{v_0,v_1,\ldots,v_r\}$ such that $d(v_0)\geqslant2,\ d(v_1)=\cdots=d(v_{r-1})=2$ and $d(v_r)=1.$ In particular, if $r=1,$ then it is a \textit{pendant edge}.  A \textit{matching} $M$ is a set of edges, no two of which have a vertex in common. In particular, $M$ is a \textit{perfect matching} if each vertex of $G$ is incident with an edge from $M$.  A subset $I$ of $V_G$ is called an \textit{independent set} if any two vertices of $I$ are not adjacent.

For two graphs $G$ and $H,$ we define $G\cup H$ their disjoint union. In addition, we use $kG$ to denote the disjoint union of $k$ copies of $G.$ The \textit{join} $G\vee H$ is the graph obtained by joining every vertex of $G$ with every vertex of $H$ with an edge. If $U\subseteq V_G,$ then we write $G[U]$ to denote the induced subgraph of $G$ with vertex set $U$ and two vertices being adjacent if and only if they are adjacent in $G.$ The \textit{complement} of a graph $G$ is a graph $\overline{G}$ with the same vertex set as $G$, in which any two distinct vertices are adjacent if and only if they are non-adjacent in $G$.

In order to formulate our main results, let $\mathcal{G}_n^\tau$\ (resp. $\mathcal{B}_n^\tau, \mathscr{T}_n^\tau$) denote the set of connected graphs (resp. bipartite graphs, trees) with order $n$ and dissociation number $\tau.$

Our first main result characterizes the $n$-vertex connected graph with dissociation number $\tau$ having the largest $A_\alpha$-index. Note that adding an edge to connect two nonadjacent vertices in a connected graph will strictly increase its $A_\alpha$-index (see, for example, \cite{V.N}). Consequently, the following theorem obviously holds.
\begin{thm}\label{thm1.1}
Let $G$ be in $\mathcal{G}_n^\tau$ having the maximum $A_\alpha$-index and let $\alpha\in [0,1)$. Then $G\cong K_{n-\tau}\vee \left(\frac{\tau}{2}K_2\right)$ if $\tau$ is even, and $G\cong K_{n-\tau}\vee \left(\frac{\tau-1}{2}K_2\cup K_1\right)$ if $\tau$ is odd.
\end{thm}

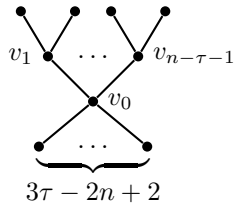
\begin{figure}[!ht]
\centering
  \begin{tikzpicture}[scale = 1.2]
  \tikzstyle{vertex}=[circle,fill=black,minimum size=0.38em,inner sep=0pt]
  \node[vertex] (G_0) at (0,0)[label=right:$v_0$]{};
  \node[vertex] (G_1) at (-0.5,0.5)[label=left:$v_1$]{};
  \node[vertex] (G_2) at (0.5,0.5)[label=right:$v_{n-\tau-1}$]{};
  \node[vertex] (G_3) at (-0.8,1){};
  \node[vertex] (G_4) at (-0.2,1){};
  \node[vertex] (G_5) at (0.2,1){};
  \node[vertex] (G_6) at (0.8,1){};
  \node[vertex] (G_7) at (-0.6,-0.5){};
  \node[vertex] (G_8) at (0.6,-0.5){};
  \draw[thick] (G_3) -- (G_1)--(G_0)--(G_8);
  \draw[thick] (G_4)--(G_1);
  \draw[thick] (G_0) -- (G_7);
  \draw[thick] (G_5)--(G_2)--(G_0);
  \draw[thick] (G_6) -- (G_2);
  \draw (0,0.5)node{$\ldots$};
  \draw (0,-0.5)node{$\ldots$};
  \node[below] at (0,-0.5){$\underbrace{\hspace{1.5cm}}$};
  \node[below] at (0,-0.8){$3\tau-2n+2$};
  \end{tikzpicture}
  \caption{Tree $S_{n,\tau}^\dag$ together with some labeled vertices.}\label{fig1b}
\end{figure}

Our second main result establishes a sharp upper bound on the $A_\alpha$-index of bipartite graphs with given order and dissociation number. The corresponding extremal graph is also characterized.
\begin{thm}\label{thm1.2}
Let $G$ be a graph in $\mathcal{B}_n^\tau.$ Assume that $\alpha\in[0,1),$ then
$$\lambda_{\alpha}(G)\leqslant \frac{1}{2}\left(\alpha n+\sqrt{\alpha^2n^2+4\tau(n-\tau)(1-2\alpha)}\right)$$
with equality if and only if $G\cong K_{\tau,n-\tau}$.
\end{thm}
Let $S^\dag_{n,\tau}$ be a tree obtained from the star $S_{n-\tau}$ by attaching exactly two pendant edges to each leaf of $S_{n-\tau}$ and attaching $3\tau-2n+2$ pendant edges to the centre of $S_{n-\tau}$ (see Figure~\ref{fig1b}). The next result determines the $n$-vertex tree with dissociation number $\tau$ having the largest $A_\alpha$-index, which reads as
\begin{thm}\label{thm1.3}
Let $T$ be a tree in $\mathscr{T}_n^\tau\ (n\geqslant 3)$ having the maximum $A_\alpha$-index with $0\leqslant \alpha < 1$. Then
$$\lambda_\alpha(T)\leqslant \theta(\alpha,n,\tau)$$
with equality if and only if $T\cong S^\dag_{n,\tau}$, where $\theta(\alpha,n,\tau)$ is the largest zero of
\begin{align}
P_\alpha(x)=&x^4+\alpha(n-2\tau-6)x^3+(8\alpha^2\tau-4\alpha^2n+4\alpha\tau-2\alpha n-2\tau+n+9\alpha^2+6\alpha-3)x^2\notag\\
&+\alpha(16\alpha n-\alpha^2n-8n-28\alpha\tau+14\tau-20\alpha+10)x+2\alpha^3n-17\alpha^2n+16\alpha n-4n\notag\\
&+24\alpha^2\tau-24\alpha\tau+6\tau-2\alpha^3+17\alpha^2-16\alpha+4.\label{eq:1.2}
\end{align}
\end{thm}
The next result characterizes the structure of all the $n$-vertex connected graphs with dissociation number $\tau\ (\tau\geqslant \left\lceil\frac{2}{3}n\right\rceil)$ having the minimum $A_\alpha$-index.
\begin{thm}\label{thm1.4}
Let $G^\dag$ be a graph in $\mathcal{G}_n^\tau$ with $\tau\geqslant \left\lceil\frac{2}{3}n\right\rceil$ having the minimum $A_\alpha$-index with $0\leqslant \alpha < 1$. Then $G^\dag$ is a tree.
\end{thm}

\begin{figure}[!ht]
\centering
  \begin{tikzpicture}[scale = 1.2]
  \tikzstyle{vertex}=[circle,fill=black,minimum size=0.38em,inner sep=0pt]
  \node[vertex] (G_1) at (0.5,0){};
  \node[vertex] (G_2) at (0,-0.5){};
  \node[vertex] (G_3) at (0,0.5){};
  \node[vertex] (G_4) at (0.85,0.5){};
  \node[vertex] (G_5) at (0.85,-0.5){};
  \node[vertex] (G_6) at (1.4,0.5){};
  \node[vertex] (G_7) at (1.4,-0.5){};
  \draw[thick] (G_3)--(G_1) -- (G_4)--(G_6);
  \draw[thick] (G_2)--(G_1) -- (G_5)--(G_7);
  \draw (0,0.08)node{$\vdots$};
  \draw (-0.3,0.06)node{$k_1$};
  \draw (1.6,0.05)node{$k_2$};
  \draw (1.1,0.05)node{$\vdots$};
  \draw (0.5,-1)node{$S_{k_1,k_2}$};
  \end{tikzpicture}
  \hspace{3em}
      \begin{tikzpicture}[scale = 1.2]
  \tikzstyle{vertex}=[circle,fill=black,minimum size=0.38em,inner sep=0pt]
  \node[vertex] (G_1) at (-1,0.5){};
  \node[vertex] (G_2) at (-0.5,0.5){};
  \node[vertex] (G_3) at (-1,-0.5){};
  \node[vertex] (G_4) at (-0.5,-0.5){};
  \node[vertex] (G_5) at (0,0){};
  \node[vertex] (G_6) at (0.4,0){};
  \node[vertex] (G_7) at (0.8,0){};
  \node[vertex] (G_8) at (1.2,0){};
  \node[vertex] (G_9) at (1.7,0.5){};
  \node[vertex] (G_10) at (2.2,0.5){};
  \node[vertex] (G_11) at (1.7,-0.5){};
  \node[vertex] (G_12) at (2.2,-0.5){};
  \draw[thick] (G_1)--(G_2) -- (G_5)--(G_6)--(G_7)--(G_8)--(G_9)--(G_10);
  \draw[thick] (G_3) -- (G_4)--(G_5);
  \draw[thick] (G_12) -- (G_11)--(G_8);
  \draw (-0.72,0.08)node{$\vdots$};
  \draw (1.95,0.08)node{$\vdots$};
  \draw (-1.23,0.05)node{$r_1$};
  \draw (2.45,0.05)node{$p_1$};
  \draw (0.7,-1)node{$T^1_{r_1,p_1}$};
  \end{tikzpicture}
\hspace{3em}
      \begin{tikzpicture}[scale = 1.2]
  \tikzstyle{vertex}=[circle,fill=black,minimum size=0.38em,inner sep=0pt]
  \node[vertex] (G_1) at (-1,0.5){};
  \node[vertex] (G_2) at (-0.5,0.5){};
  \node[vertex] (G_3) at (-1,-0.5){};
  \node[vertex] (G_4) at (-0.5,-0.5){};
  \node[vertex] (G_5) at (0,0){};
  \node[vertex] (G_6) at (0.4,0){};
  \node[vertex] (G_7) at (0.8,0){};
  \node[vertex] (G_8) at (1.2,0){};
  \node[vertex] (G_9) at (1.7,0.5){};
  \node[vertex] (G_10) at (2.2,0.5){};
  \node[vertex] (G_11) at (1.7,-0.5){};
  \node[vertex] (G_12) at (2.2,-0.5){};
  \node[vertex] (G_13) at (0,0.5){};
  \draw[thick] (G_1)--(G_2) -- (G_5)--(G_6)--(G_7)--(G_8)--(G_9)--(G_10);
  \draw[thick] (G_3) -- (G_4)--(G_5);
  \draw[thick] (G_12) -- (G_11)--(G_8);
  \draw[thick] (G_13) -- (G_5);
  \draw (-0.72,0.08)node{$\vdots$};
  \draw (1.95,0.08)node{$\vdots$};
  \draw (-1.23,0.05)node{$r_2$};
  \draw (2.45,0.05)node{$p_2$};
  \draw (0.7,-1)node{$T^2_{r_2,p_2}$};
  \end{tikzpicture}
  \caption{Trees $S_{k_1,k_2}, T^1_{r_1,p_1}$ and $ T^2_{r_2,p_2}$}\label{fig01}
\end{figure}
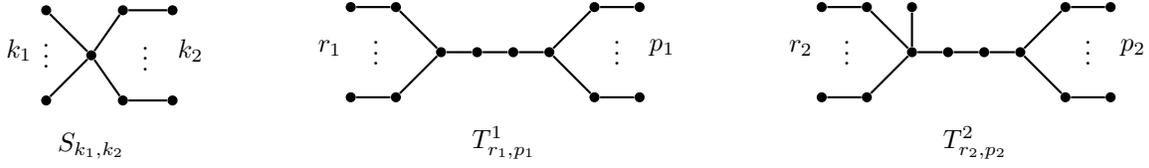

Denote by $S_{k_1,k_2}$ the tree obtained from $S_{k_1+1}$ by attaching $k_2$ pendant paths of length two to the centre of $S_{k_1+1}.$
Let $T^1_{r_1,p_1}$ be  the tree obtained from $P_4$ by attaching $r_1$ and $p_1$ pendant paths of length two to the two leaves of $P_4$, respectively. Then let $T^2_{r_2,p_2}$ be the tree obtained from $T^1_{r_2,p_2}$ by attaching one pendant edge to the vertex of degree $r_2+1$ in $T^1_{r_2,p_2}$. The trees $S_{k_1,k_2}, T^1_{r_1,p_1}$ and $ T^2_{r_2,p_2}$ are depicted in Figure~\ref{fig01}. 

Our last main result determines all the connected $n$-vertex graphs with dissociation number $\tau\in\{2,\left\lceil\frac{2}{3}n\right\rceil,n-1,n-2\}$ having the minimum $A_\alpha$-index.
\begin{thm}\label{thm1.5}
Let $G^*$ be a graph in $\mathcal{G}_n^\tau$ having the minimum $A_\alpha$-index with $0\leqslant\alpha<1$.
\begin{wst}
\item[{\rm (i)}] If $\tau=2,$ then $G^*\cong K_n-M$, where $M$ is a maximum matching of $K_n.$
\item[{\rm (ii)}] If $\tau=\left\lceil\frac{2}{3}n\right\rceil,$ then $G^*\cong P_n.$
\item[{\rm (iii)}] If $n\geqslant4$ and $\tau=n-1,$ then $G^*\cong S_{0,\frac{n-1}{2}}$ if $n$ is odd and $G^*\cong S_{1,\frac{n-2}{2}}$ if $n$ is even.
\item[{\rm (iv)}] If $n\geqslant6$  and $\tau=n-2,$ then $G^*\cong T^1_{\left\lceil\frac{n-4}{4}\right\rceil,\left\lfloor\frac{n-4}{4}\right\rfloor}$ if $n$ is even and $G^*\cong T^2_{\left\lfloor\frac{n-5}{4}\right\rfloor,\left\lceil\frac{n-5}{4}\right\rceil}$ if $n$ is odd.
\end{wst}
\end{thm}
The remainder of this paper is organized as follows. In Section 2, we give some essential definitions and some necessary preliminaries. In Section 3, we give the proofs of Theorems \ref{thm1.2} and \ref{thm1.3}. In Section 4, we give the proofs for Theorems \ref{thm1.4} and \ref{thm1.5}. Some concluding remarks are given in the last section.
{\section{\normalsize Preliminaries}\setcounter{equation}{0}

In this section, we give some necessary results, which will be used to prove our main results.
\begin{lem}[\cite{V.N}]\label{lem2.1}
If $G$ is a connected graph and $H$ is a proper subgraph of $G$, then $\lambda_\alpha(H)<\lambda_\alpha(G)$ for $\alpha\in[0,1).$
\end{lem}

\begin{lem}[\cite{C.Y}]\label{lem2.2}
If $T$ is a tree of order $n$ and $\alpha\in[0,1],$ then $\lambda_\alpha(T)\leqslant \frac{1}{2}\left(\alpha n+\sqrt{\alpha^2 n^2+4(n-1)(1-2\alpha)}\right)$ with equality if and only if $T\cong S_n.$
\end{lem}
\begin{lem}[\cite{C.Y}]\label{lem2.3}
If $G$ is a connected graph of order $n$ and $\alpha\in[0,1],$ then $\lambda_\alpha(G)\geqslant \lambda_\alpha(P_n)$ with equality if and only if $G\cong P_n.$
\end{lem}
For $V_1\subseteq V_G$ and $E_1\subseteq E_G,$ we write $G-V_1$ and $G-E_1$ for the graphs obtained from $G$ by deleting the vertices in $V_1$ and their incident edges and the edges in $E_1$, respectively. For convenience, we may use $G-v$ and $G-uv$ to denote $G-\{v\}$ and $G-\{uv\},$ respectively. Similarly, $G+uv$ is obtained from $G$ by adding the edge $uv\notin E_G.$

By the well-known Perron-Frobenius theorem, if $G$ is connected, we know that the multiplicity of $\lambda_\alpha(G)$ is one and there exists a unit positive eigenvector, say ${\bf x}=(x_1,\ldots,x_n)^T,$ corresponding to it, where $x_i$ is the component of ${\bf x}$ at $v_i\in V_G$ \cite{V.N}. As usual we call ${\bf x}$  the \textit{Perron vector} of $A_\alpha(G).$
\begin{lem}[\cite{HZ}]\label{lem2.4}
Let $G$ be a connected graph with $u,v\in V_G$ and let $\alpha\in[0,1).$ Assume that $\bf x$ is the Perron vector of $A_\alpha(G)$ such that $x_u\geqslant x_v.$ If $\{v_1,v_2,\ldots,v_s\}\in N_G(v)\setminus N_G(u)\ (1\leqslant s\leqslant d_v)$ and $G'=G-\{vv_1,\ldots,vv_s\}+\{uv_1,\ldots,uv_s\},$ then $\lambda_\alpha(G')>\lambda_\alpha(G).$
\end{lem}
\begin{lem}[\cite{CS}]\label{lem2.5}
Let $G$ be a connected graph with $|E_G|\geqslant1$ and $u\in V_G.$ Denote by $G_{s,t}$ the graph obtained from $G$ by attaching two pendant paths of length $s$ and $t$ to $u$, respectively. If $s\geqslant t\geqslant 1$ and $0\leqslant\alpha<1$, then $\lambda_\alpha(G_{s,t})>\lambda_\alpha(G_{s+1,t-1}).$
\end{lem}
Let $P=v_1v_2\ldots v_k$ with $k\geqslant2$ be a path in a connected graph. Then we call $P$ an \textit{internal path}, if $\min\{d(v_1),d(v_k)\}\geqslant3$ and $d(v_2)=\cdots=d(v_{k-1})=2.$
The \textit{subdivision operation} for an edge $uv\in E_G$ is adding a new vertex $w$ and substituting $uv$ by a path $uwv,$ and we denote the resultant graph by $G_w.$
\begin{lem}[\cite{FLZ}]\label{lem2.6}
Let $G$ be a connected graph with $\alpha\in[0,1)$ and $uv$ be an edge on an internal path of $G$. Let $G_w$ be the graph obtained from $G$ by the subdivision operation for the edge $uv$. Then $\lambda_\alpha(G_w)<\lambda_\alpha(G).$
\end{lem}

Let $R$ be a real matrix, whose rows and columns are indexed by $V=\{1,2,\ldots,n\}$. Assume that $\pi=\{V_1,V_2,\ldots,V_t\}$ is a partition of $V$. Then $R$ can be partitioned based on $\pi$ as
\begin{align*}
R=\begin{pmatrix}
R_{11} & \cdots & R_{1t} \\
\vdots &  \ddots & \vdots \\
R_{t1} & \cdots & R_{tt}
\end{pmatrix},
\end{align*}
where $R_{ij}$ denotes the submatrix of $R$, indexed by the rows and columns of $V_i$ and $V_j$, respectively. Let $r_{ij}$ be the average row sum of $R_{ij}$ for $1\leqslant i, j\leqslant t$. Usually, the $t\times t$ matrix $R^{\pi}=(r_{ij})$ is called the \textit{quotient matrix} of $R$. Moreover, if the row sum of $R_{ij}$ is constant for $1\leqslant i,j\leqslant t,$ then we call $\pi$ an \textit{equitable partition}.
\begin{lem}[\cite{YYSX}]\label{lem2.7}
Let $R$ be a real matrix with an equitable partition $\pi$, and let $R^\pi$ be the corresponding quotient matrix. Then each eigenvalue of $R^\pi$ is an eigenvalue of $R.$ Furthermore, if $R$ is nonnegative, then the spectral radii of $R$ and $R^\pi$ are equal.
\end{lem}

Recall that $\mathcal{P}_G$ denotes the set of all pendant vertices of $G$ and $\mathcal{Q}_G'$ denotes the set of all quasi-pendant vertices of degree 2 of $G.$ The following lemma reveals that there exists a maximum dissociation set of $G$ containing all vertices in $\mathcal{P}_G\bigcup\mathcal{Q}_G',$ which plays an essential role in the proofs of our main results.
\begin{lem}[\cite{HLZ}]\label{lem2.8}
Let $G$ be a graph with order $n\geqslant5.$ Then there exists a maximum dissociation set $S$ such that $\mathcal{P}_G\bigcup \mathcal{Q}_G'\subseteq S.$
\end{lem}
A maximum dissociation set is said to be \textit{good} if it contains all vertices in $\mathcal{P}_G\bigcup\mathcal{Q}_G'.$ According to Lemma~\ref{lem2.8}, we know that such a set always exists.
\section{\normalsize Proofs of Theorems \ref{thm1.2} and \ref{thm1.3}}\setcounter{equation}{0}
In this section, we give the proofs for Theorems \ref{thm1.2} and \ref{thm1.3}. The former identifies the unique $n$-vertex bipartite graph with dissociation number $\tau$ having the largest $A_\alpha$-index, and the later characterizes the tree with given order $n$ and dissociation number $\tau$ having the maximum $A_\alpha$-index.
\begin{proof}[\bf Proof of Theorem \ref{thm1.2}]
Choose $G=(X,Y)\in\mathcal{B}_n^\tau$ such that its $A_\alpha$-index is as large as possible. We first assume that $G$ is connected. Without loss of generality, we assume that $|X|\geqslant|Y|.$ Let $S$ be a maximum dissociation set of $G$. Since $X$ is a dissociation set of $G$, we have $\tau=|S|\geqslant|X|.$ Clearly, our result is true for $n\leqslant 2.$ So in what follows, we consider $n\geqslant 3$.

If $\tau=|X|,$ by Lemma \ref{lem2.1}, we obtain that $G\cong K_{\tau,n-\tau}.$ In the following, we may assume that $\tau>|X|.$ Therefore, $S\cap Y\neq\emptyset.$  Let $X_1=X\cap S,\ Y_2=Y\cap S$ and $X_2=X\setminus X_1,\ Y_1=Y\setminus Y_2.$ For convenience, assume $|X_1|=a,\ |Y_1|=b,\ |X_2|=c,$ and $|Y_2|=d.$ Since $|X_1|+|Y_2|=|S|>|X|\geqslant|Y|$, we have $d>c,\ a>b.$ According to the choice of $G$ and by Lemma~\ref{lem2.1}, one sees that both $G[(X_1, Y_1)]$ and $G[(X_2, Y)]$ are complete, and $G[(X_1, Y_2)]$ contains independent edges as many as possible. We proceed by considering the following three cases according to the values of $a$ and $d.$

{\bf Case 1.}\ $a=d.$ In this case, $E_{G[X_1\cup Y_2]}$ is a perfect matching of $G[X_1\cup Y_2]$. It is easy to see that $\pi_1:=X_1\cup X_2\cup Y_1\cup Y_2$ is an equitable partition of $V_G,$ and the corresponding quotient matrix can be written as follows:
$$
(A_\alpha)^{\pi_1}:=\left(
                 \begin{array}{cccc}
                   \alpha(b+1) & 0 & (1-\alpha)b & 1-\alpha  \\
                   0 & \alpha(a+b) & (1-\alpha)b & (1-\alpha)a \\
                   (1-\alpha)a & (1-\alpha)c & \alpha(a+c) & 0 \\
                   1-\alpha & (1-\alpha)c & 0 & \alpha(c+1)  \\
                 \end{array}
               \right).
$$
Let $P_\alpha^1(x):=\det\left(x I_4-(A_\alpha)^{\pi_1}\right)$  be the characteristic polynomial of $(A_\alpha)^{\pi_1}.$ Note that $(A_\alpha)^{\pi_1}$ is nonnegative and irreducible. Together with Lemma~\ref{lem2.7} and the Perron-Frobenius theorem, we know that $\lambda_\alpha(G)$ coincides with the largest root of $P_\alpha^1(x)=0.$

Since $n\geqslant3$ and $G$ is connected, we have $a=d>\max\{b,c\}\geqslant1.$ Recall that $|X|\geqslant|Y|$ and $a=d,$ we can deduce that $c\geqslant b.$ Let $G'\cong K_{2a,b+c}.$ It is easy to see that $G'\in\mathcal{B}_n^\tau$. By a simple calculation, we have
$$
\lambda_\alpha(G')=\frac{1}{2}\left(\alpha(2a+b+c)+\sqrt{\alpha^2(2a+b+c)^2+8a(b+c)(1-2\alpha)}\right).
$$
Utilizing calculations by Matlab, we obtain that $\min P_\alpha^1(x)\approx 2.598$ for $x\geqslant\lambda_\alpha(G')$ (see the Appendix), i.e. $P_\alpha^1(x)>0$ when $x\geqslant\lambda_\alpha(G').$ Recall that $\lambda_\alpha(G)$ is the largest zero of $P_\alpha^1(x).$ That is to say, $\lambda_\alpha(G)<\lambda_\alpha(G'),$ which contradicts the choice of $G.$


{\bf Case 2.}\ $a<d.$ In this case, there exists a set $Y_2'\subseteq Y_2$ such that $|Y_2'|=a$ and $E_{G[X_1\cup Y_2']}$ is a perfect matching of $G[X_1\cup Y_2']$. Let $Y_2''=Y_2\setminus Y_2',$ then $Y_2'' \neq\emptyset.$ Thus it is easy to see that $\pi_2:=X_1\cup X_2\cup Y_1\cup Y_2'\cup Y_2''$ is an equitable partition of $V_G,$ and the corresponding quotient matrix can be written as follows:
$$
(A_\alpha)^{\pi_2}:=\left(
                 \begin{array}{ccccc}
                   \alpha(b+1) & 0 & (1-\alpha)b & 1-\alpha & 0  \\
                   0 & \alpha(b+d) & (1-\alpha)b & (1-\alpha)a & (1-\alpha)(d-a) \\
                   (1-\alpha)a & (1-\alpha)c & \alpha(a+c) & 0 & 0 \\
                   1-\alpha & (1-\alpha)c & 0 & \alpha(c+1) & 0  \\
                   0 & (1-\alpha)c & 0 & 0 & \alpha c  \\
                 \end{array}
               \right).
$$
Let $P_\alpha^2(x):=\det\left(x I_4-(A_\alpha)^{\pi_2}\right)$  be the characteristic polynomial of $(A_\alpha)^{\pi_2}.$ Note that $(A_\alpha)^{\pi_2}$ is nonnegative and irreducible. Together with Lemma~\ref{lem2.7} and the Perron-Frobenius theorem, we know that $\lambda_\alpha(G)$ coincides with the largest root of $P_\alpha^2(x)=0.$

It is easy to see that $c\geqslant1$ (since otherwise $|X|=|X_1|=a<d=|Y_2|\leqslant|Y|,$ which contradicts $|X|\geqslant|Y|$). Thus we have $d\geqslant c+1\geqslant2.$ Let $G''\cong K_{a+d,b+c}.$ It is easy to see that $G''\in\mathcal{B}_n^\tau$. By a simple calculation, we have
$$
\lambda_\alpha(G'')=\frac{1}{2}\left(\alpha(a+b+c+d)+\sqrt{\alpha^2(a+b+c+d)^2+4(a+d)(b+c)(1-2\alpha)}\right).
$$
With the same idea as the Appendix of calculations using Matlab, we obtain that $\min P_\alpha^2(x)=1$ for $x\geqslant\lambda_\alpha(G'')$, that is to say, $P_\alpha^2(G)>0$ when $x\geqslant\lambda_\alpha(G'').$  Recall that $\lambda_\alpha(G)$ is the largest zero of $P_\alpha^2(x).$ That is to say, $\lambda_\alpha(G)<\lambda_\alpha(G''),$ which contradicts the choice of $G.$

{\bf Case 3.}\ $a>d.$ We can also derive a contradiction in a similar way as Case 2, whose procedure is omitted here.

By Cases 1-3, we obtain $G\cong K_{\tau,n-\tau}.$ A simple calculation gives us
$$\lambda_{\alpha}(G)=\frac{1}{2}\left(\alpha n+\sqrt{\alpha^2n^2+4\tau(n-\tau)(1-2\alpha)}\right).$$

In what follows, we show $G$ is connected. Otherwise suppose that $G$ is disconnected.
Let $G=\bigcup_{i=1}^sG_i,$ where $G_i$ is a connected bipartite graph in $\mathcal{B}_{n_i}^{\tau_i}$ and $s\geqslant2$. Thus we have
$$
\sum_{i=1}^s(n_i-\tau_i)=\sum_{i=1}^sn_i-\sum_{i=1}^s\tau_i=n-\tau.
$$
Therefore, $n_i<n, \tau_i<\tau$ and $n_i-\tau_i<n-\tau.$ Consequently, by Lemma~\ref{lem2.1},
$$
\lambda_\alpha(G)=\max\left\{\lambda_\alpha(G_i)\,|\,1\leqslant i\leqslant s\right\}
                 \leqslant\max\left\{\lambda_\alpha(K_{\tau_i,n_i-\tau_i})\,|\,1\leqslant i\leqslant s\right\}
                 <\lambda_\alpha(K_{\tau,n-\tau}),
$$
a contradiction to the choice of $G.$

This completes the proof.
\end{proof}
Let $\mathscr{T}^1_{n,\tau}$ be a set of all the $n$-vertex trees obtained from $S_{n-\tau+1}$ by attaching exactly one pendant edge to the centre of $S_{n-\tau+1}$, and attaching at least two pendent edges to each leaf of $S_{n-\tau+1}$ such that the total number of leaves of the resultant tree is $\tau-1$, where $\tau\notin\left\{\frac{2}{3}n,\frac{2n+1}{3}\right\}$.
Let $\mathscr{T}^2_{n,\tau}$ be a set of all the $n$-vertex trees obtained from $S_{n-\tau+1}$ by attaching at least two pendent edges to every leaf of $S_{n-\tau+1}$ such that the total number of leaves of the resultant tree is $\tau-1$, where $\tau\neq \frac{2}{3}n$.
Let $\mathscr{T}^3_{n,\tau}$ be a set of all the $n$-vertex trees obtained from $S_{n-\tau}$ by attaching at least two pendant edges to each vertex of $S_{n-\tau}$ such that the total number of leaves of the resultant tree is $\tau.$ In Figure~\ref{fig02}, one sees a tree in $\mathscr{T}^1_{n,\tau}, \mathscr{T}^2_{n,\tau}$ and $\mathscr{T}^3_{n,\tau},$  respectively. Let $\mathcal{S}_{n,\tau}=\mathscr{T}^1_{n,\tau}\bigcup\mathscr{T}^2_{n,\tau}\bigcup\mathscr{T}^3_{n,\tau}.$ In view of Lemma~\ref{lem2.8}, we have $\mathcal{S}_{n,\tau}\subseteq\mathscr{T}_n^\tau.$ Sun and Li \cite{LS} showed that $\tau(F)\geqslant\left\lceil\frac{2}{3}n\right\rceil$ for each forest $F$ with order $n.$ Therefore, for every tree $T$ in $\mathcal{S}_{n,\tau},$ $\tau(T)\geqslant\left\lceil\frac{2}{3}n\right\rceil,$ then $3\tau-2n+2\geqslant2,$ which implies that $S^\dag_{n,\tau}\in\mathscr{T}^3_{n,\tau}.$

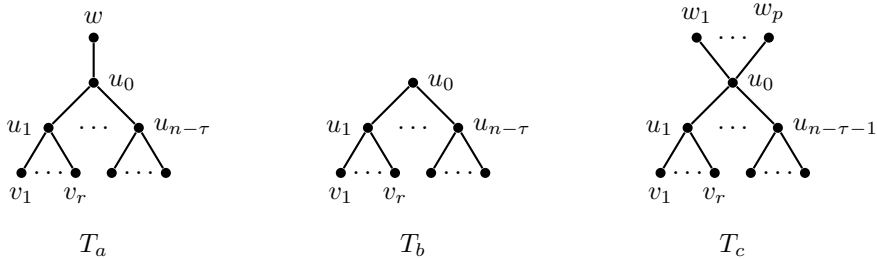
\begin{figure}[!ht]
\centering
\begin{tikzpicture}[scale = 1.2]
  \tikzstyle{vertex}=[circle,fill=black,minimum size=0.38em,inner sep=0pt]
  \node[vertex] (G_0) at (0,0)[label=right:$u_0$]{};
  \node[vertex] (G_1) at (-0.5,-0.5)[label=left:$u_1$]{};
  \node[vertex] (G_2) at (0.5,-0.5)[label=right:$u_{n-\tau}$]{};
  \node[vertex] (G_3) at (-0.8,-1)[label=below:$v_1$]{};
  \node[vertex] (G_4) at (-0.2,-1)[label=below:$v_r$]{};
  \node[vertex] (G_5) at (0.2,-1){};
  \node[vertex] (G_6) at (0.8,-1){};
  \node[vertex] (G_7) at (0,0.5)[label=above:$w$]{};
  \draw[thick] (G_3) -- (G_1)--(G_0)--(G_7);
  \draw[thick] (G_4)--(G_1);
  \draw[thick] (G_5)--(G_2)--(G_0);
  \draw[thick] (G_6) -- (G_2);
  \draw (0,-0.5)node{$\ldots$};
  \draw (0.5,-1)node{$\ldots$};
  \draw (-0.5,-1)node{$\ldots$};
  \draw (0,-1.8)node{$T_a$};
  \end{tikzpicture}
   \hspace{3em}
\begin{tikzpicture}[scale = 1.2]
  \tikzstyle{vertex}=[circle,fill=black,minimum size=0.38em,inner sep=0pt]
  \node[vertex] (G_0) at (0,0)[label=right:$u_0$]{};
  \node[vertex] (G_1) at (-0.5,-0.5)[label=left:$u_1$]{};
  \node[vertex] (G_2) at (0.5,-0.5)[label=right:$u_{n-\tau}$]{};
  \node[vertex] (G_3) at (-0.8,-1)[label=below:$v_1$]{};
  \node[vertex] (G_4) at (-0.2,-1)[label=below:$v_r$]{};
  \node[vertex] (G_5) at (0.2,-1){};
  \node[vertex] (G_6) at (0.8,-1){};
  \draw[thick] (G_3) -- (G_1)--(G_0);
  \draw[thick] (G_4)--(G_1);
  \draw[thick] (G_5)--(G_2)--(G_0);
  \draw[thick] (G_6) -- (G_2);
  \draw (0,-0.5)node{$\ldots$};
  \draw (0.5,-1)node{$\ldots$};
  \draw (-0.5,-1)node{$\ldots$};
  \draw (0,-1.8)node{$T_b$};
  \end{tikzpicture}
  \hspace{3em}
  \begin{tikzpicture}[scale = 1.2]
  \tikzstyle{vertex}=[circle,fill=black,minimum size=0.38em,inner sep=0pt]
  \node[vertex] (G_0) at (0,0)[label=right:$u_0$]{};
  \node[vertex] (G_1) at (-0.5,-0.5)[label=left:$u_1$]{};
  \node[vertex] (G_2) at (0.5,-0.5)[label=right:$u_{n-\tau-1}$]{};
  \node[vertex] (G_3) at (-0.8,-1)[label=below:$v_1$]{};
  \node[vertex] (G_4) at (-0.2,-1)[label=below:$v_r$]{};
  \node[vertex] (G_5) at (0.2,-1){};
  \node[vertex] (G_6) at (0.8,-1){};
  \node[vertex] (G_7) at (-0.4,0.5)[label=above:$w_1$]{};
  \node[vertex] (G_8) at (0.4,0.5)[label=above:$w_p$]{};
  \draw[thick] (G_3) -- (G_1)--(G_0)--(G_7);
  \draw[thick] (G_4)--(G_1);
  \draw[thick] (G_5)--(G_2)--(G_0);
  \draw[thick] (G_6) -- (G_2);
  \draw[thick] (G_8) -- (G_0);
  \draw (0,-0.5)node{$\ldots$};
  \draw (0.5,-1)node{$\ldots$};
  \draw (-0.5,-1)node{$\ldots$};
  \draw (0,0.5)node{$\ldots$};
  \draw (0,-1.8)node{$T_c$};
  \end{tikzpicture}
  \caption{Trees $T_a\in \mathscr{T}^1_{n,\tau}, \ T_b\in\mathscr{T}^2_{n,\tau}$ and $T_c\in\mathscr{T}^3_{n,\tau}$ together with some labeled vertices.}\label{fig02}
\end{figure}

In order to prove Theorem~\ref{thm1.3}, we need the following lemma.
\begin{lem}\label{lem3.1}
For $0\leqslant\alpha<1$, the $A_\alpha$-index of $S^\dag_{n,\tau}$ is equal to the largest root of $P_\alpha(x)=0$, where $P_\alpha(x)$ is defined in \eqref{eq:1.2}.
\end{lem}
\begin{proof}
Recall that $S^\dag_{n,\tau}$ is a tree obtained from the star $S_{n-\tau}$ by attaching exactly two pendant edges to each leaf of $S_{n-\tau}$ and attaching $3\tau-2n+2$ pendant edges to the centre of $S_{n-\tau}$ (see Figure~\ref{fig1b}). Let $v_0$ be the centre of $S_{n-\tau}$, $V_1$ be the set of $2(n-\tau-1)$ pendant vertices adjacent to the leaves of $S_{n-\tau}$, and let $V_2$ be the set of $3\tau-2n+2$ pendant vertices adjacent to $v_0$. Hence $\pi:=\{v_0\}\cup (V_{S_{n-\tau}}\setminus\{v_0\}) \cup V_1 \cup V_2$ is an equitable partition of $V_{S^\dag_{n,\tau}}$. Then the corresponding quotient matrix can be written as follows:
$$
(A_\alpha)^{\pi}:=\left(
                 \begin{array}{cccc}
                   \alpha(2\tau-n+1) & (1-\alpha)(n-\tau-1) & 0 & (1-\alpha)(3\tau-2n+2)  \\
                   1-\alpha & 3\alpha & 2(1-\alpha) & 0 \\
                   0 & 1-\alpha & \alpha & 0 \\
                   1-\alpha & 0 & 0 & \alpha  \\
                 \end{array}
               \right).
$$
The characteristic polynomial of $(A_\alpha)^{\pi}$ is $P_\alpha(x),$ where $P_\alpha(x)$ is defined in \eqref{eq:1.2}. Note that $(A_\alpha)^{\pi}$ is nonnegative and irreducible. Together with Lemma~\ref{lem2.7} and the Perron-Frobenius theorem, our result holds immediately.
\end{proof}

\begin{proof}[\bf Proof of Theorem \ref{thm1.3}]
Choose $T$ among $\mathscr{T}_n^\tau$ such that its $A_\alpha$-index is as large as possible. Let ${\bf x}$ be the Perron vector of $A_\alpha(T).$ In view of Lemma~\ref{lem2.2}, we obtain that $T\cong S_n\cong S^\dag_{n,n-1}$ if $\tau=n-1$ or $3\leqslant n\leqslant4$. So, in what follows, we assume that $\tau\leqslant n-2$ and $n\geqslant5.$ We are to characterize the structure of $T$ through the following three claims.
\begin{claim}\label{c1}
$\mathcal{Q}_T'=\emptyset$.
\end{claim}
\begin{proof}[\bf Proof of Claim \ref{c1}]
Suppose that $\mathcal{Q}_T'\neq\emptyset.$ Let $v\in\mathcal{Q}_T'$ and $N_T(v)=\{u,w\},$ where $u$ is a leaf. Denote by $T_1=T-vu+uw.$ Then by Lemma~\ref{lem2.5}, we have $\lambda_\alpha(T_1)>\lambda_\alpha(T).$ However, by Lemma~\ref{lem2.8}, $T_1\in\mathscr{T}_n^\tau,$ a contradiction to the choice of $T.$
Thus $\mathcal{Q}_T'=\emptyset.$ This completes the proof of Claim~\ref{c1}.
\end{proof}
\begin{claim}\label{c2}
$T\in\mathcal{S}_{n,\tau}$.
\end{claim}
\begin{proof}[\bf Proof of Claim \ref{c2}]
Let $T'=T-\mathcal{P}_T.$ By Claim~\ref{c1}, one sees that $T\in\mathcal{S}_{n,\tau}$ is equivalent to $|\mathcal{Q}_{T'}|=1.$ Suppose that $|\mathcal{Q}_{T'}|\geqslant2.$ Choose two vertices $u,v\in\mathcal{Q}_{T'}$ such that the distance between $u$ and $v$ in $T'$ is as large as possible. Let $u'$ and $v'$ be the neighbors of $u$ and $v$ lying on the path connecting $u$ and $v$, respectively. Then we have $\left(N_T(u)\setminus\{u'\}\right)\cup \left(N_T(v)\setminus\{v'\}\right)\subseteq\mathcal{P}_T\cup\mathcal{Q}_T,$ $\left(N_T(u)\setminus\{u'\}\right)\cap\mathcal{Q}_T\neq\emptyset$ and $\left(N_T(v)\setminus\{v'\}\right)\cap\mathcal{Q}_T\neq\emptyset.$

Assume that $\left(N_T(u)\setminus\{u'\}\right)\cap\mathcal{Q}_T=\{u_1,u_2,\ldots,u_s\}$ and $\left(N_T(v)\setminus\{v'\}\right)\cap\mathcal{Q}_T=\{v_1,v_2,\ldots,v_t\}.$ Thus we have $s,t\geqslant1.$ By Claim~\ref{c1}, we know that there exist at least two leaves adjacent to each vertex in $\{u_1,\ldots,u_s,v_1,\ldots,v_t\}.$ Let $S$ be a good dissociation set of $T.$ Then $u_i\notin S$ and $v_j\notin S$ for $1\leqslant i\leqslant s,\ 1\leqslant j\leqslant t.$ Let
\begin{equation*}
    T_2=
    \left\{
    \begin{array}{ll}
        T-\{vv_j\,|\,1\leqslant j\leqslant t\}
        +\{uv_j\,|\,1\leqslant j\leqslant t\},& \textrm{if $x_u>x_v$,}\\[5pt]
        T-\{uu_i\,|\,1\leqslant i\leqslant s\}
        +\{vu_i\,|\,1\leqslant i\leqslant s\},& \textrm{if $x_u\leqslant x_v.$}
    \end{array}
    \right.
\end{equation*}
In view of Lemmas~\ref{lem2.4} and \ref{lem2.8}, we have $\lambda_\alpha(T_2)>\lambda_\alpha(T)$ and $T_2\in\mathscr{T}_n^\tau,$ a contradiction. Thus we have $|\mathcal{Q}_{T'}|=1,$ i.e. $T\in\mathcal{S}_{n,\tau}$.
\end{proof}
\begin{claim}\label{c3}
$T\in\mathscr{T}^3_{n,\tau}$.
\end{claim}
\begin{proof}[\bf Proof of Claim \ref{c3}]
Suppose that $T\notin\mathscr{T}^3_{n,\tau}.$ Then by Claim~\ref{c2}, we know that $T\in\mathscr{T}^1_{n,\tau}\bigcup\mathscr{T}^2_{n,\tau}.$ For the star $S_{n-\tau+1},$ assume that $u_0$ is the centre of $S_{n-\tau+1},$ and $u_1,\ldots,u_{n-\tau}$ are leaves of $S_{n-\tau+1}$ (see Figure~\ref{fig02}).  We proceed by considering the following two cases.

{\bf Case 1.}\ $T\in\mathscr{T}^1_{n,\tau}.$ Assume that $N(u_1)\setminus\{u_0\}=\{v_1,\ldots,v_r\}$ and $w$ is just the leaf adjacent to $u_0$ in $T$. Let
\begin{equation*}
    T_3=
    \left\{
    \begin{array}{ll}
        T-\{u_1v_i\,|\,1\leqslant i\leqslant r\}
        +\{u_0v_i\,|\,1\leqslant i\leqslant r\},& \textrm{if $x_{u_0}>x_{u_1}$,}\\[5pt]
        T-u_0w+u_1w,& \textrm{if $x_{u_0}\leqslant x_{u_1}.$}
    \end{array}
    \right.
\end{equation*}
By Lemma~\ref{lem2.8}, we obtain $T_3\in\mathscr{T}^2_{n,\tau}\bigcup\mathscr{T}^3_{n,\tau}.$ However, in view of Lemma~\ref{lem2.4}, $\lambda_\alpha(T_3)>\lambda_\alpha(T)$, which contradicts the choice of $T.$

{\bf Case 2.}\ $T\in\mathscr{T}^2_{n,\tau}.$ Similarly, assume that $N(u_1)\setminus\{u_0\}=\{v_1,\ldots,v_r\}.$ Let
\begin{equation*}
    T_4=
    \left\{
    \begin{array}{ll}
        T-\{u_1v_i\,|\,1\leqslant i\leqslant r\}
        +\{u_0v_i\,|\,1\leqslant i\leqslant r\},& \textrm{if $x_{u_0}>x_{u_1}$,}\\[5pt]
        T-\{u_0u_j\,|\,2\leqslant j\leqslant n-\tau\}
        +\{u_1u_j\,|\,2\leqslant j\leqslant n-\tau\},& \textrm{if $x_{u_0}\leqslant x_{u_1}.$}
    \end{array}
    \right.
\end{equation*}
By Lemma~\ref{lem2.8}, we get that $T_4\in\mathscr{T}^3_{n,\tau}.$ And once again, by Lemma~\ref{lem2.4}, $\lambda_\alpha(T_4)>\lambda_\alpha(T)$, a contradiction.
\end{proof}
We come back to show Theorem~\ref{thm1.3}.

By Claim~\ref{c3}, we have $T\in\mathscr{T}^3_{n,\tau}$. Similarly, for the star $S_{n-\tau}$, assume that $V_{S_{n-\tau}}=\{u_0,u_1,\ldots,u_{n-\tau-1}\}$, where $u_0$ is the unique non-pendant vertex. Suppose that $T\ncong S^\dag_{n,\tau}.$ Then there exists a vertex $u_i,$ where $i\in\{1,2,\ldots,n-\tau-1\},$ such that $d(u_i)\geqslant4.$ Without loss of generality, assume that $d(u_1)\geqslant4.$
Let $N(u_0)\cap \mathcal{P}_T=\{w_1,\ldots,w_p\}$ and $N(u_1)\setminus\{u_0\}=\{v_1,\ldots,v_r\}$. Then $p\geqslant2$ and $r\geqslant3.$

If $\tau=n-2,$ then $T$ is a tree obtained from $S_{p+1}$ and $S_{r+1}$ by adding an edge to connect the centres of theirs. Since $T\ncong S^\dag_{n,\tau},$ we have $p\geqslant3.$  Let
\begin{equation*}
    T_5=
    \left\{
    \begin{array}{ll}
        T-\{u_1v_i\,|\,3\leqslant i\leqslant r\}
        +\{u_0v_i\,|\,3\leqslant i\leqslant r\},& \textrm{if $x_{u_0}>x_{u_1}$,}\\[7pt]
        T-\{u_0w_j\,|\,3\leqslant j\leqslant p\}
        +\{u_1w_j\,|\,3\leqslant j\leqslant p\},& \textrm{if $x_{u_0}\leqslant x_{u_1}.$}
    \end{array}
    \right.
\end{equation*}
Then $T_5\cong S^\dag_{n,n-2}$. By Lemma~\ref{lem2.4}, we have $\lambda_\alpha(T_5)>\lambda_\alpha(T)$, a contradiction.

If $\tau\leqslant n-3,$ then let
\begin{eqnarray*}
    T_6=
    \left\{
    \begin{array}{ll}
        T-\{u_1v_i\,|\,3\leqslant i\leqslant r\}+\{u_0v_i\,|\,3\leqslant i\leqslant r\},& \textrm{if $x_{u_0}>x_{u_1},$}\\[5pt]
        T-\{u_0u_k\,|\,2\leqslant k\leqslant n-\tau-1\}+
        \{u_1u_k\,|\,2\leqslant k\leqslant n-\tau-1\},& \textrm{if $x_{u_0}\leqslant x_{u_1}$ and $p=2,$}\\[5pt]
        T-\{u_0u_k\,|\,2\leqslant k\leqslant n-\tau-1\}-\{u_0w_j\,|\,3\leqslant j\leqslant p\}\\[5pt]
         +\{u_1u_k\,|\,2\leqslant k\leqslant n-\tau-1\}+\{u_1w_j\,|\,3\leqslant j\leqslant p\},& \textrm{if $x_{u_0}\leqslant x_{u_1}$ and $p\geqslant3.$}
    \end{array}
    \right.
\end{eqnarray*}
By Lemma~\ref{lem2.8}, $\tau(T_6)=\tau(T).$ In view of Lemma~\ref{lem2.4}, one has $\lambda_\alpha(T_6)>\lambda_\alpha(T)$, a contradiction.
Therefore, $T\cong S^\dag_{n,\tau}.$ By Lemma~\ref{lem3.1}, our result holds immediately.

This completes the proof.
\end{proof}
{\section{\normalsize Proofs of Theorems \ref{thm1.4} and \ref{thm1.5}}\setcounter{equation}{0}
In this section, we present the proofs for Theorems~\ref{thm1.4} and \ref{thm1.5}. The former shows that a connected graph with given order $n$ and dissociation number $\tau\ (\tau\geqslant\left\lceil\frac{2}{3}n\right\rceil)$ having the minimum $A_\alpha$-index is a tree, and the latter characterizes all the connected graphs with order $n$ and dissociation number $\tau\in\{2,\left\lceil\frac{2}{3}n\right\rceil,n-1,n-2\}$  having the minimum $A_\alpha$-index.

Let $Y_1$ (resp. $Y_2$) be a tree obtained from $P_{n-2}$ (resp. $P_{n-6}$) by attaching exactly two pendant edges (resp. two pendant paths of length three) to one leaf of $P_{n-2}$ (resp. $P_{n-6}$). Let $Y_3$ be a tree obtained from $P_{n-4}$ by attaching one pendant edge and one pendant path of length three to one leaf of $P_{n-4}.$ Recall that $G_w$ is a graph constructed from $G$ by a subdivision operation for an edge of $G$. The \textit{triple subdivision operation} for an edge $uv\in E_G$ is adding three new vertices $x,y,z$ and substituting $uv$ by a path $uxyzv,$ and we denote the resultant graph by $G_{xyz}.$ We call $\dot{G}$ a \textit{subdivision transformation graph} of $G$ if $\dot{G}$ is obtained from $G$ by a (triple) subdivision operation for an edge on an internal path of $G,$ and deleting one or three other vertices of $G$ such that $|V_{\dot{G}}|=|V_G|.$ A subdivision transformation graph $\dot{G}$ is said to be \textit{optimal} if $\dot{G}$ is connected and $\tau(\dot{G})\in\{\tau(G)-1,\tau(G)\}.$

In order to show Theorem~\ref{thm1.4}, we need the following key lemma. It presents the relationship between $\tau(G_w)$ (resp. $\tau(G_{xyz})$) and $\tau(G).$
\begin{lem}[\cite{HLZ}]\label{lem4.1}
Let $G$ be a connected graph with $uv\in E_G$ and let $G_w$ and $G_{xyz}$ be the graph obtained from $G$ by the subdivision operation and triple subdivision operation for $uv$, respectively. Then
\begin{wst}
\item[{\rm (i)}] $\tau(G_w)\in\{\tau(G), \tau(G)+1\};$
\item[{\rm (ii)}] $\tau(G_{xyz})=\tau(G)+2.$
\end{wst}
\end{lem}
Let $T$ be a tree, we call $v\in V_T$ a \textit{branching vertex} if $d_T(v)\geqslant3$. The next lemma was obtained previously in \cite{HLZ}. We include its proof here.
\begin{lem}[\cite{HLZ}]\label{lem4.2}
Let $T$ be a tree with at least two branching vertices. Then there exists an optimal subdivision transformation graph of $T$.
\end{lem}
\begin{proof}
We use $T_x$ and $T_{xyz}$ to denote the trees obtained from $T$ by the subdivision operation and triple subdivision operation for an edge on an internal path of $T,$ respectively. Choose a diameter path $P_d=u_1u_2u_3u_4\cdots u_d$ of $T$ such that $d_T(u_2)$ is as large as possible. Then $N_T(u_2)\setminus\{u_3\}\subseteq \mathcal{P}_T$ and $N_T(u_3)\setminus\{u_4\}\subseteq\mathcal{P}_T\cup\mathcal{Q}_T.$

If $d_T(u_2)=d_T(u_3)=2,$ then by Lemma~\ref{lem2.8}, we obtain
\begin{align}\label{eq4.1}
\tau(T_{xyz}-u_1-u_2-u_3)=\tau(T_{xyz})-2.
\end{align}
By Lemma~\ref{lem4.1}, we have $\tau(T_{xyz})=\tau(T)+2.$ Combining this with \eqref{eq4.1} yields $\tau(T_{xyz}-u_1-u_2-u_3)=\tau(T).$ Therefore, $T_{xyz}-u_1-u_2-u_3$ is an optimal subdivision transformation graph of $T.$

If $d_T(u_2)=2$ and $d_T(u_3)\geqslant3,$ then according to the choice of $P_d,$ we know that $N_T(u_3)\setminus\{u_4\}\subseteq\mathcal{P}_T\cup\mathcal{Q}_T'.$ Hence, in view of Lemmas~\ref{lem2.8} and \ref{lem4.1}, we obtain that $T_x-u_1$ is a tree with $\tau(T_x-u_1)=\tau(T_x)-1\in\{\tau(T),\tau(T)-1\}.$ This implies that $T_x-u_1$ is an optimal subdivision transformation graph of $T.$

If $d_T(u_2)=3,$ then by Lemmas~\ref{lem2.8} and \ref{lem4.1}, we have $T_{xyz}-u_1-u_2-v$ is a tree with $\tau(T_{xyz}-u_1-u_2-v)=\tau(T_{xyz})-2=\tau(T),$ where $v$ is the only pendant vertex adjacent to $u_2$ other than $u_1.$ This indicates that $T_{xyz}-u_1-u_2-v$ is an optimal subdivision transformation graph of $T.$

If $d_T(u_2)\geqslant4,$ then in view of Lemmas~\ref{lem2.8} and \ref{lem4.1}, we have $T_x-u_1$ is a tree with $\tau(T_x-u_1)=\tau(T_x)-1\in\{\tau(T),\tau(T)-1\}.$ Hence $T_x-u_1$ is an optimal subdivision transformation graph of $T.$

This completes the proof.
\end{proof}
Now we are ready to show Theorem~\ref{thm1.4}.
\begin{proof}[\bf Proof of Theorem \ref{thm1.4}]
If $\tau(G^\dag)=\left\lceil\frac{2}{3}n\right\rceil,$ then by Lemma~\ref{lem2.3}, we have $G^\dag\cong P_n$ since $\tau(P_n)=\left\lceil\frac{2}{3}n\right\rceil,$ and our result holds obviously. So, in what follows, we assume that $\tau(G^\dag)>\left\lceil\frac{2}{3}n\right\rceil.$

Suppose to the contrary that $G^\dag$ is not a tree. If $G^\dag$ is a subgraph of $(nK_2)\vee K_1,$ then $\tau(G^\dag)=n-1$ and there exists a cycle $v_0v_1v_2v_0$ such that $d(v_1)=d(v_2)=2.$ Let $G^\dag_1=G^\dag-v_1v_2.$ By Lemmas~\ref{lem2.1} and \ref{lem2.8}, we have $\lambda_\alpha(G^\dag_1)<\lambda_\alpha(G^\dag)$ and  $G^\dag_1\in\mathcal{G}_n^{n-1},$ which contradicts the choice of $G^\dag$. This means that there exists a spanning tree, say $T^\dag,$ of $G^\dag$ such that $T^\dag\ncong S_{k_1,k_2}$ (see Figure~\ref{fig01}) for all $k_1+2k_2=n-1.$

It is clear that $\tau(T^\dag)\geqslant\tau(G^\dag).$ By Lemma~\ref{lem2.1}, we have
\begin{align}\label{Eq:4.2}
\lambda_\alpha(T^\dag)<\lambda_\alpha(G^\dag).
\end{align}

If $\tau(T^\dag)=\tau(G^\dag),$ then we can derive a contradiction by the choice of $G^\dag.$ Hence,
\begin{align}\label{eq:4.1}
\tau(T^\dag)>\tau(G^\dag)\geqslant\left\lceil\frac{2}{3}n\right\rceil+1.
\end{align}

Since $\tau(T^\dag)>\left\lceil\frac{2}{3}n\right\rceil=\tau(P_n),$ we have $T^\dag\ncong P_n,$ which implies that there exists at least one branching vertex in $T^\dag.$

If $T^\dag$ has at least two branching vertices, then by Lemma~\ref{lem4.2}, there exists an optimal subdivision transformation graph, say $\dot{T}_1,$ of $T^\dag.$ Denote by $\dot{T}_0=T^\dag.$ In view of Lemma~\ref{lem2.6} and \eqref{Eq:4.2}, we obtain $\lambda_\alpha(\dot{T}_1)<\lambda_\alpha(\dot{T}_0)<\lambda_\alpha(G^\dag).$
We can repeat the above transformation to get an $n$-vertex tree sequence
$$
\dot{T}_0, \dot{T}_1,\ldots,\dot{T}_i,\ldots,\dot{T}_s
$$
such that $\dot{T}_s$ is a tree having exactly one branching vertex with $\tau(\dot{T}_s)=\left\lceil\frac{2}{3}n\right\rceil$ and $\tau(\dot{T}_i)\in\{\tau(\dot{T}_{i-1})-1, \tau(\dot{T}_{i-1})\}$ and
$\lambda_\alpha(\dot{T}_i)<\lambda_\alpha(\dot{T}_{i-1})$ for $1\leqslant i\leqslant s$. According to the proof of Lemma~\ref{lem4.2}, we know that such $\dot{T}_s$ must exist, for example, $Y_1,Y_2,Y_3$.
By \eqref{eq:4.1}, $\tau(\dot{T}_0)\geqslant\lceil\frac{2}{3}n\rceil+2=\tau(\dot{T}_s)+2$, which implies that $s\geqslant2$. Since $\tau(\dot{T}_0)>\tau(G^\dag)>\left\lceil\frac{2}{3}n\right\rceil
=\tau(\dot{T}_s),$  there exists $j\in \{1,2,\ldots,s-1\}$ such that $\tau(\dot{T}_j)=\tau(G^\dag).$ However,
$\lambda_\alpha(\dot{T}_j)<\lambda_\alpha(\dot{T}_0)<\lambda_\alpha(G^\dag),$ which contradicts the choice of $G^\dag.$
\begin{figure}[!ht]
\centering
  \begin{tikzpicture}[scale = 1.8]
  \tikzstyle{vertex}=[circle,fill=black,minimum size=0.38em,inner sep=0pt]
  \node[vertex] (G_0) at (0,0)[label=above:$u_0$]{};
  \node[vertex] (G_1) at (-0.7,-0.3)[label=left:$u^1_1$]{};
  \node[vertex] (G_2) at (-0.7,-0.6)[label=left:$u^1_2$]{};
  \node[vertex] (G_3) at (-0.7,-0.93){};
  \node[vertex] (G_4) at (-0.7,-1.23)[label=left:$u^1_{n_1}$]{};
  \node[vertex] (G_5) at (-0.3,-0.3)[label=right:$u^2_1$]{};
  \node[vertex] (G_6) at (-0.3,-0.6)[label=right:$u^2_2$]{};
  \node[vertex] (G_7) at (-0.3,-0.93){};
  \node[vertex] (G_8) at (-0.3,-1.23)[label=right:$u^2_{n_2}$]{};
  \node[vertex] (G_9) at (0.5,-0.3)[label=right:$u^t_1$]{};
  \node[vertex] (G_10) at (0.5,-0.6)[label=right:$u^t_2$]{};
  \node[vertex] (G_11) at (0.5,-0.93){};
  \node[vertex] (G_12) at (0.5,-1.23)[label=right:$u^t_{n_t}$]{};
  \draw[thick] (G_0)--(G_1)--(G_2);
  \draw[thick] (G_0)--(G_5)--(G_6);
  \draw[thick] (G_0)--(G_9)--(G_10);
  \draw[thick] (G_4)--(G_3);
  \draw[thick] (G_8) -- (G_7);
  \draw[thick] (G_12)--(G_11);
  \draw (-0.7,-0.7)node{$\vdots$};
  \draw (-0.3,-0.7)node{$\vdots$};
  \draw (0.5,-0.7)node{$\vdots$};
  \draw (0.1,-0.45)node{$\ldots$};
  \draw (0.1,-1.05)node{$\ldots$};
  \draw (-0.1,-1.7)node{$T^\dag$};
  \end{tikzpicture}
		\hspace{3em}
 \begin{tikzpicture}[scale = 1.8]
  \tikzstyle{vertex}=[circle,fill=black,minimum size=0.38em,inner sep=0pt]
  \node[vertex] (G_1) at (0,-0.4)[label=below:$u_0$]{};
  \node[vertex] (G_2) at (-0.3,0)[label=above:$u_1$]{};
  \node[vertex] (G_3) at (0.3,0)[label=above:$u_r$]{};
  \node[vertex] (G_4) at (-0.3,-0.8)[label=left:$u_{r+1}$]{};
  \node[vertex] (G_5) at (-0.3,-1.2)[label=left:$u_{r+t+1}$]{};
  \node[vertex] (G_6) at (0.3,-0.8)[label=right:$u_{r+t}$]{};
  \node[vertex] (G_7) at (0.3,-1.2)[label=right:$u_{r+2t}$]{};
  \node[vertex] (G_8) at (0.5,-0.4){};
  \node[vertex] (G_9) at (1,-0.4)[label=below:$u_{n-3}$]{};
  \node[vertex] (G_10) at (1.5,-0.4)[label=below:$u_{n-2}$]{};
  \node[vertex] (G_11) at (2,-0.4)[label=below:$u_{n-1}$]{};
  \draw[thick] (G_3)--(G_1) -- (G_4)--(G_5);
  \draw[thick] (G_2)--(G_1) -- (G_6)--(G_7);
  \draw[thick] (G_1)--(G_8);
  \draw[thick] (G_9)--(G_10) -- (G_11);
  \draw (0,0)node{$\ldots$};
  \draw (0,-1)node{$\ldots$};
  \draw (0.75,-0.4)node{$\ldots$};
  \draw (0.85,-1.7)node{$T^\dag_2\cong W_{r,t}$};
  \end{tikzpicture}
  \caption{Trees $T^\dag$ and $T^\dag_2$.}\label{fig03}
\end{figure}
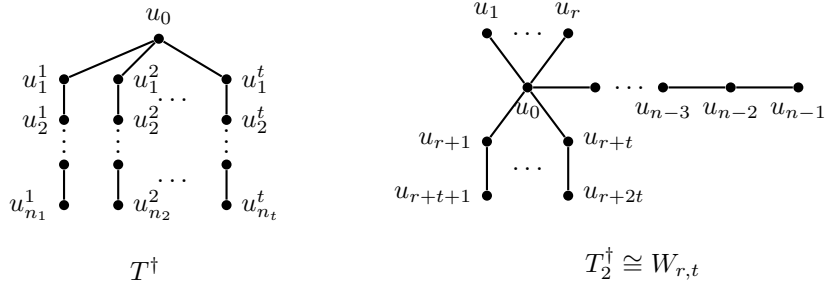

In what follows, we assume that $T^\dag$ has exact one branching vertex. Let $u_0$ be the branching vertex of $T^\dag.$ Then $T^\dag$ can be obtained from $t$ paths $P_{n_1+1},P_{n_2+1},\ldots,P_{n_t+1}$ by sharing a common vertex $u_0,$ where $t\geqslant3.$ Let $u_0,u^i_1,\ldots,u^i_{n_i}$ denote the vertices of $P_{n_i+1}$ for $1\leqslant i\leqslant t$ (see Figure~\ref{fig03}). Without loss of generality, we can assume that $n_1\geqslant n_2\geqslant\cdots\geqslant n_t\geqslant1.$ Since $T^\dag\ncong S_{k_1,k_2}$ for all $k_1+2k_2=n-1,$ we have $n_1\geqslant3.$ Let
\begin{equation*}
    T^\dag_1=
    \left\{
    \begin{array}{ll}
        T^\dag,& \textrm{if $n_2\in\{1,2\},$}\\[3pt]
        T^\dag-u_0u^2_1+u^1_{n_1}u^2_1,& \textrm{if $n_2\equiv 0\pmod 3,$}\\[3pt]
        T^\dag-u^2_1u^2_2+u^1_{n_1}u^2_2,& \textrm{if $n_2\equiv 1\pmod 3$ and $n_2\geqslant4,$}\\[3pt]
        T^\dag-u^2_2u^2_3+u^1_{n_1}u^2_3,& \textrm{if $n_2\equiv 2\pmod 3$ and $n_2\geqslant5.$}
    \end{array}
    \right.
\end{equation*}
In view of Lemmas~\ref{lem2.5} and \ref{lem2.8}, we get that $\tau(T^\dag_1)=\tau(T^\dag)$ and $\lambda_\alpha(T^\dag_1)\leqslant \lambda_\alpha(T^\dag).$ Repeating the above procedure for $T^\dag_1,$ we finally obtain an $n$-vertex tree $T^\dag_2\cong W_{r,t},$ such that
\begin{align}\label{eq:4.2}
\lambda_\alpha(T^\dag_2)\leqslant \lambda_\alpha(T^\dag),\ \ \tau(T^\dag_2)=\tau(T^\dag),
\end{align}
where $W_{r,t}$ is an $n$-vertex tree obtained from a star $S_{r+1}$ by attaching $t$ pendant paths of length 2 and one pendant path of length $n-r-2t-1$ to the centre of $S_{r+1}$ (see Figure~\ref{fig03}).

If $(r,t)=(3,0)$ or $r+t\leqslant2$, then we can obtain
\begin{equation*}
    \tau(T^\dag_2)=
    \left\{
    \begin{array}{ll}
        \left\lceil\frac{2}{3}n\right\rceil,& \textrm{if $(r,t)\in\{(0,0),(0,1),(1,0),(2,0)\},$}\\[6pt]
        \left\lceil\frac{2n-1}{3}\right\rceil+1,& \textrm{if $(r,t)=(0,2),$}\\[5pt]
        \left\lceil\frac{2n-2}{3}\right\rceil+1,& \textrm{if $(r,t)\in\{(1,1),(3,0)\}.$}
    \end{array}
    \right.
\end{equation*}
By \eqref{eq:4.1} and \eqref{eq:4.2}, we have $\tau(T^\dag_2)=\tau(T^\dag)\geqslant\left\lceil\frac{2}{3}n\right\rceil+2,$ a contradiction.

If $t=0$ and $r\geqslant4,$ then by Lemmas~\ref{lem2.5} and \ref{lem2.8}, we have
\begin{align}\label{eq:4.05}
\lambda_\alpha(W_{r-2,t+1})< \lambda_\alpha(W_{r,t}),\ \ \tau(W_{r-2,t+1})=\tau(W_{r,t}).
\end{align}
In view of \eqref{eq:4.05}, we assume that $t\geqslant1$ and $r+t\geqslant3$ in what follows. The labels of the vertices of $T^\dag_2$ are shown in Figure~\ref{fig03}. Let
\begin{align*}
    T^\dag_3=
    \left\{
    \begin{array}{ll}
        T^\dag_2-\{u_0u_{r+j}\,|\,1\leqslant j\leqslant t-2\}+\{u_{n-3}u_{r+j}\,|\,1\leqslant j\leqslant t-2\},& \textrm{if $r=0,$ $t\geqslant3,$}\\[5pt]
        T^\dag_2-\{u_0u_{r+j}\,|\,1\leqslant j\leqslant t-1\}+\{u_{n-3}u_{r+j}\,|\,1\leqslant j\leqslant t-1\},& \textrm{if $r=1,$ $t\geqslant2,$}\\[5pt]
        T^\dag_2-\{u_0u_i\,|\,1\leqslant i\leqslant r-1\}+\{u_{n-3}u_i\,|\,1\leqslant i\leqslant r-1\},& \textrm{if $r\geqslant2,$ $t=1,$}\\[5pt]
        T^\dag_2-\{u_0u_{r+j}\,|\,1\leqslant j\leqslant t-1\}-\{u_0u_i\,|\,1\leqslant i\leqslant r-1\}\\[5pt]
+\{u_{n-3}u_{r+j}\,|\,1\leqslant j\leqslant t-1\}+\{u_{n-3}u_i\,|\,1\leqslant i\leqslant r-1\},& \textrm{if $r,t\geqslant2.$}
    \end{array}
    \right.
\end{align*}
In view of Lemma~\ref{lem2.8}, we have
\begin{align}\label{eq:4.3}
\tau(T^\dag_3)=\tau(T^\dag_2).
\end{align}
Let ${\bf y}=(y_0,y_1,\ldots,y_{n-1})^T$ be the Perron vector of $T^\dag_3$, where $y_i$ corresponds to $u_i\ (0\leqslant i\leqslant n-1).$  If $y_0<y_{n-3},$ let
\begin{align*}
    T^\dag_4=
    \left\{
    \begin{array}{ll}
        T^\dag_3-u_0u_{r+t-1}+u_{n-3}u_{r+t-1},& \textrm{if $r=0,$}\\[5pt]
        T^\dag_3-u_0u_r+u_{n-3}u_r,& \textrm{if $r>0.$}
    \end{array}
    \right.
\end{align*}
It is obvious that $T^\dag_4\cong T^\dag_2.$ By Lemma~\ref{lem2.4}, we have $\lambda_\alpha(T^\dag_3)<\lambda_\alpha(T^\dag_4)=\lambda_\alpha(T^\dag_2).$
If $y_0\geqslant y_{n-3},$ then we can undo the step from $T_2^\dag$ to $T_3^\dag.$
By Lemma~\ref{lem2.4}, we have $\lambda_\alpha(T^\dag_3)<\lambda_\alpha(T^\dag_2)$. Consequently, both cases yield $\lambda_\alpha(T^\dag_3)<\lambda_\alpha(T^\dag_2).$ Together with \eqref{Eq:4.2}, \eqref{eq:4.1}, \eqref{eq:4.2} and \eqref{eq:4.3}, we have
\begin{align}\label{eq:4.4}
\lambda_\alpha(T^\dag_3)< \lambda_\alpha(G^\dag),\ \ \tau(T^\dag_3)>\tau(G^\dag)\geqslant\left\lceil\frac{2}{3}n\right\rceil+1.
\end{align}
Note that $T^\dag_3$ has two branching vertices. We obtain, similarly as in the case ``$T^\dag$ has at least two branching vertices," that there exists an $n$-vertex tree sequence
$$
\ddot{T}_0, \ddot{T}_1,\ldots,\ddot{T}_i,\ldots,\ddot{T}_p
$$
such that $\ddot{T}_0\cong T^\dag_3,$ $\ddot{T}_p$ is a tree having exactly one branching vertex with $\tau(\dot{T}_s)=\left\lceil\frac{2}{3}n\right\rceil$ and $\tau(\ddot{T}_i)\in\{\tau(\ddot{T}_{i-1})-1, \tau(\ddot{T}_{i-1})\}, \lambda_\alpha(\ddot{T}_i)<\lambda_\alpha(\ddot{T}_{i-1})$ for $1\leqslant i\leqslant p$. Note that $p\geqslant2.$ By \eqref{eq:4.4}, we have $\tau(\ddot{T}_0)>\tau(G^\dag)>\tau(\ddot{T}_p).$ This implies that there exists $j\in\{1,2,\ldots,p-1\}$ such that $\tau(\ddot{T}_j)=\tau(G^\dag).$ However, we know that $\lambda_\alpha(\ddot{T}_j)<\lambda_\alpha(\ddot{T}_0)<\lambda_\alpha(G^\dag),$ which contradicts the choice of $G^\dag.$

This completes the proof.
\end{proof}

For an $n$-vertex graph $G$, denote by $P_\alpha(G,x)=\det(x I_n-A_\alpha(G))$ the characteristic polynomial of $A_\alpha(G)$, where $I_n$ is the identity
matrix of order $n$. The following result is obvious.
\begin{lem}\label{lem2.9}
Let $G_1$ and $G_2$ be two connected graphs. If $P_\alpha(G_2,x)>P_\alpha(G_1,x)$ for $x\geqslant\lambda_\alpha(G_1)$, then $\lambda_\alpha(G_2)<\lambda_\alpha(G_1)$.
\end{lem}

Recall that $T^1_{r_1,p_1}$ is a tree obtained from $P_4$ by attaching $r_1$ and $p_1$ pendant paths of length two to the two leaves of $P_4$, respectively, whereas $T^2_{r_2,p_2}$ is a tree obtained from $T^1_{r_2,p_2}$ by attaching one pendant edge to the vertex of degree $r_2+1$ in $T^1_{r_2,p_2}$ (see Figure~\ref{fig01}). In order to show Theorem~\ref{thm1.5}, we need the following lemma.
\begin{lem}\label{lem4.4}
If $r_1\geqslant p_1\geqslant1,$ then $\lambda_{\alpha}(T_{r_1,p_1}^1)<\lambda_{\alpha}(T_{r_1+1,p_1-1}^1).$
\end{lem}
\begin{proof}
For simplicity, let $G^*:=T_{r_1,p_1}^1$ and $G':=T_{r_1+1,p_1-1}^1.$ If $p_1=1,$ then the result follows immediately from Lemma~\ref{lem2.5}. So in what follows, we may assume that $r_1\geqslant p_1\geqslant2.$ We can easily obtain an equitable partition of $V_{G^*}$ (resp. $V_{G'}$), with four parts each of which contains only one vertex, two parts each of which contains $r_1$ (resp. $r_1+1$) vertices and two parts each of which contains $p_1$ (resp. $p_1-1$) vertices. Therefore, the corresponding quotient matrices can be written as follows:
$$
(A_\alpha)^{\pi_3}:=
\left(
                 \begin{array}{cccccccc}
                   \alpha & 1-\alpha & 0 & 0 & 0 & 0 & 0 & 0  \\
                   1-\alpha & 2\alpha & 1-\alpha & 0 & 0 & 0 & 0 & 0 \\
                   0 & (1-\alpha)r_1 & \alpha(r_1+1) & 1-\alpha & 0 & 0 & 0 & 0 \\
                   0 & 0 & 1-\alpha & 2\alpha & 1-\alpha & 0 & 0 & 0  \\
                   0 & 0 & 0 & 1-\alpha & 2\alpha & 1-\alpha & 0 & 0  \\
                   0 & 0 & 0 & 0 & 1-\alpha & \alpha(p_1+1) & (1-\alpha)p_1 & 0  \\
                   0 & 0 & 0 & 0 & 0 & 1-\alpha & 2\alpha & 1-\alpha \\
                   0 & 0 & 0 & 0 & 0 & 0 & 1-\alpha & \alpha \\
                 \end{array}
               \right),
$$

$$
(A_\alpha)^{\pi_4}:
=\left(
                 \begin{array}{cccccccc}
                   \alpha & 1-\alpha & 0 & 0 & 0 & 0 & 0 & 0  \\
                   1-\alpha & 2\alpha & 1-\alpha & 0 & 0 & 0 & 0 & 0 \\
                   0 & (1-\alpha)(r_1+1) & \alpha(r_1+2) & 1-\alpha & 0 & 0 & 0 & 0 \\
                   0 & 0 & 1-\alpha & 2\alpha & 1-\alpha & 0 & 0 & 0  \\
                   0 & 0 & 0 & 1-\alpha & 2\alpha & 1-\alpha & 0 & 0  \\
                   0 & 0 & 0 & 0 & 1-\alpha & \alpha p_1 & (1-\alpha)(p_1-1) & 0  \\
                   0 & 0 & 0 & 0 & 0 & 1-\alpha & 2\alpha & 1-\alpha \\
                   0 & 0 & 0 & 0 & 0 & 0 & 1-\alpha & \alpha \\
                 \end{array}
               \right).
$$
Let $P_\alpha^3(x):=\det\left(x I_8-(A_\alpha)^{\pi_3}\right)$ (resp. $P_\alpha^4(x):=\det\left(x I_8-(A_\alpha)^{\pi_4}\right)$)  be the characteristic polynomial of $(A_\alpha)^{\pi_3}$ (resp. $(A_\alpha)^{\pi_4}$). Note that $(A_\alpha)^{\pi_3}$ (resp. $(A_\alpha)^{\pi_4}$) is nonnegative and
irreducible. Together with Lemma~\ref{lem2.7} and the Perron-Frobenius theorem, we know that $\lambda_\alpha(G^*)$ (resp. $\lambda_\alpha(G')$) coincides with the largest zero of $P_\alpha^3(x)$ (resp. $P_\alpha^4(x)$).

Note that $G'$ contains $S_{r_1+3}$ as its proper subgraph, then by Lemma~\ref{lem2.1}, we have
\begin{align}\label{eq4.8}
\lambda_\alpha(G')>\lambda_\alpha(S_{r_1+3}),
\end{align}
and a simple calculation gives us
$$\lambda_\alpha(S_{r_1+3})=\frac{1}{2}\left(\alpha(r_1+3)+\sqrt{\alpha^2(r_1+3)^2+4(r_1+2)(1-2\alpha)}\right).$$
With the same idea as the Appendix of calculations using Matlab, we obtain that $\min\left(P_\alpha^3(x)-P_\alpha^4(x)\right)\approx0.6527$ for $x\geqslant\lambda_\alpha(S_{r_1+3})$, that is to say, $P_\alpha^3(x)-P_\alpha^4(x)>0$ when $x\geqslant\lambda_\alpha(S_{r_1+3})$. Combining this with \eqref{eq4.8} yields $P_\alpha^3(x)-P_\alpha^4(x)>0$ when $x\geqslant\lambda_\alpha(G').$ Then by Lemma~\ref{lem2.9}, we have $\lambda_\alpha(G^*)<\lambda_\alpha(G'),$ as desired.
\end{proof}

Let $T^3_{r_3,p_3}$ (resp. $T^4_{r_4,p_4}$) be the tree obtained from $S_4$ (resp. $P_2$) by attaching $r_3$ (resp. $r_4$) and $p_3$ (resp. $p_4$) pendant paths of length two to two leaves of $S_4$ (resp. $P_2$), respectively. Let $T^5_{r_5,p_5}$ be the tree obtained from $P_3$ by attaching one pendant edge and $r_5$ pendant paths of length two to one leaf of $P_3$ and attaching $p_5$ pendant paths of length two to the other leaf of $P_3.$ Let $T^6_{r_6,p_6}$ be the tree obtained from $S_{1,2}$ by attaching $r_6$ and $p_6$ pendant paths of length two to the two quasi-pendant vertices with degree 2 of $S_{1,2}$, respectively.
Let $T^7_{r_7,p_7}$ (resp. $T^8_{r_8,p_8}$) be the tree obtained from $P_4$ (resp. $P_6$) by attaching $r_7$ (resp. $r_8$) and $p_7$ (resp. $p_8$) pendant paths of length two to two quasi-pendant vertices of $P_4$ (resp. $P_6$), respectively.  $T^i_{r_i,p_i}$ for $i\in\{3,\ldots,8\}$ are depicted in Figure~\ref{fig04}.

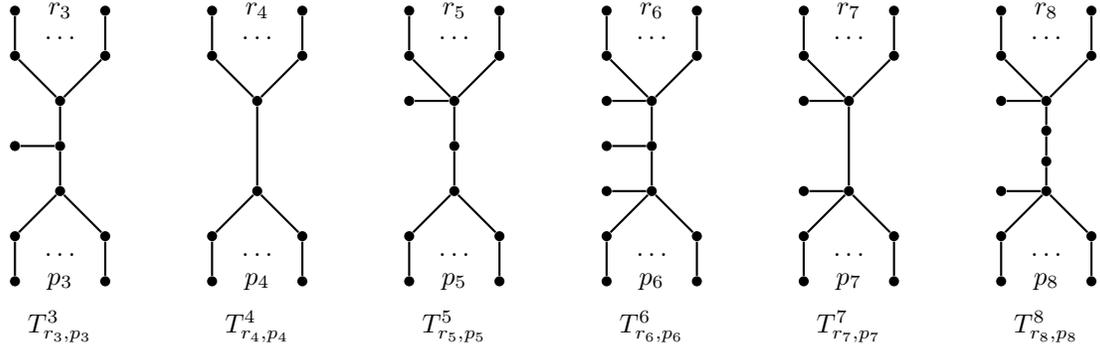
\begin{figure}[!ht]
\centering
      \begin{tikzpicture}[scale = 1.2]
  \tikzstyle{vertex}=[circle,fill=black,minimum size=0.38em,inner sep=0pt]
  \node[vertex] (G_1) at (0,0){};
  \node[vertex] (G_2) at (-0.5,0){};
  \node[vertex] (G_3) at (0,0.5){};
  \node[vertex] (G_4) at (-0.5,1){};
  \node[vertex] (G_5) at (-0.5,1.5){};
  \node[vertex] (G_6) at (0.5,1.5){};
  \node[vertex] (G_7) at (0.5,1){};
  \node[vertex] (G_8) at (0,-0.5){};
  \node[vertex] (G_9) at (-0.5,-1){};
  \node[vertex] (G_10) at (-0.5,-1.5){};
  \node[vertex] (G_11) at (0.5,-1){};
  \node[vertex] (G_12) at (0.5,-1.5){};
  \draw[thick] (G_5)--(G_4) -- (G_3)--(G_1)--(G_8)--(G_9)--(G_10);
  \draw[thick] (G_6)--(G_7)--(G_3);
  \draw[thick] (G_12)--(G_11)--(G_8);
  \draw[thick] (G_2)--(G_1);
  \draw (0,1.2)node{$\ldots$};
  \draw (0,-1.2)node{$\ldots$};
  \draw (0,1.5)node{$r_3$};
  \draw (0,-1.5)node{$p_3$};
  \draw (0,-2)node{$T^3_{r_3,p_3}$};
  \end{tikzpicture}
     \hspace{3em}
  \begin{tikzpicture}[scale = 1.2]
  \tikzstyle{vertex}=[circle,fill=black,minimum size=0.38em,inner sep=0pt]
  \node[vertex] (G_3) at (0,0.5){};
  \node[vertex] (G_4) at (-0.5,1){};
  \node[vertex] (G_5) at (-0.5,1.5){};
  \node[vertex] (G_6) at (0.5,1.5){};
  \node[vertex] (G_7) at (0.5,1){};
  \node[vertex] (G_8) at (0,-0.5){};
  \node[vertex] (G_9) at (-0.5,-1){};
  \node[vertex] (G_10) at (-0.5,-1.5){};
  \node[vertex] (G_11) at (0.5,-1){};
  \node[vertex] (G_12) at (0.5,-1.5){};
  \draw[thick] (G_5)--(G_4) -- (G_3)--(G_8)--(G_9)--(G_10);
  \draw[thick] (G_6)--(G_7)--(G_3);
  \draw[thick] (G_12)--(G_11)--(G_8);
  \draw (0,1.2)node{$\ldots$};
  \draw (0,-1.2)node{$\ldots$};
  \draw (0,1.5)node{$r_4$};
  \draw (0,-1.5)node{$p_4$};
  \draw (0,-2)node{$T^4_{r_4,p_4}$};
  \end{tikzpicture}
     \hspace{3em}
           \begin{tikzpicture}[scale = 1.2]
  \tikzstyle{vertex}=[circle,fill=black,minimum size=0.38em,inner sep=0pt]
  \node[vertex] (G_1) at (0,0){};
  \node[vertex] (G_2) at (-0.5,0.5){};
  \node[vertex] (G_3) at (0,0.5){};
  \node[vertex] (G_4) at (-0.5,1){};
  \node[vertex] (G_5) at (-0.5,1.5){};
  \node[vertex] (G_6) at (0.5,1.5){};
  \node[vertex] (G_7) at (0.5,1){};
  \node[vertex] (G_8) at (0,-0.5){};
  \node[vertex] (G_9) at (-0.5,-1){};
  \node[vertex] (G_10) at (-0.5,-1.5){};
  \node[vertex] (G_11) at (0.5,-1){};
  \node[vertex] (G_12) at (0.5,-1.5){};
  \draw[thick] (G_5)--(G_4) -- (G_3)--(G_1)--(G_8)--(G_9)--(G_10);
  \draw[thick] (G_6)--(G_7)--(G_3);
  \draw[thick] (G_12)--(G_11)--(G_8);
  \draw[thick] (G_2)--(G_3);
  \draw (0,1.2)node{$\ldots$};
  \draw (0,-1.2)node{$\ldots$};
  \draw (0,1.5)node{$r_5$};
  \draw (0,-1.5)node{$p_5$};
  \draw (0,-2)node{$T^5_{r_5,p_5}$};
  \end{tikzpicture}
  \hspace{3em}
        \begin{tikzpicture}[scale = 1.2]
  \tikzstyle{vertex}=[circle,fill=black,minimum size=0.38em,inner sep=0pt]
  \node[vertex] (G_1) at (0,0){};
  \node[vertex] (G_2) at (-0.5,0){};
  \node[vertex] (G_3) at (0,0.5){};
  \node[vertex] (G_4) at (-0.5,1){};
  \node[vertex] (G_5) at (-0.5,1.5){};
  \node[vertex] (G_6) at (0.5,1.5){};
  \node[vertex] (G_7) at (0.5,1){};
  \node[vertex] (G_8) at (0,-0.5){};
  \node[vertex] (G_9) at (-0.5,-1){};
  \node[vertex] (G_10) at (-0.5,-1.5){};
  \node[vertex] (G_11) at (0.5,-1){};
  \node[vertex] (G_12) at (0.5,-1.5){};
  \node[vertex] (G_13) at (-0.5,0.5){};
  \node[vertex] (G_14) at (-0.5,-0.5){};
  \draw[thick] (G_5)--(G_4) -- (G_3)--(G_1)--(G_8)--(G_9)--(G_10);
  \draw[thick] (G_6)--(G_7)--(G_3);
  \draw[thick] (G_12)--(G_11)--(G_8);
  \draw[thick] (G_2)--(G_1);
  \draw[thick] (G_13)--(G_3);
  \draw[thick] (G_14)--(G_8);
  \draw (0,1.2)node{$\ldots$};
  \draw (0,-1.2)node{$\ldots$};
  \draw (0,1.5)node{$r_6$};
  \draw (0,-1.5)node{$p_6$};
  \draw (0,-2)node{$T^6_{r_6,p_6}$};
  \end{tikzpicture}
       \hspace{3em}
  \begin{tikzpicture}[scale = 1.2]
  \tikzstyle{vertex}=[circle,fill=black,minimum size=0.38em,inner sep=0pt]
  \node[vertex] (G_1) at (-0.5,0.5){};
  \node[vertex] (G_2) at (-0.5,-0.5){};
  \node[vertex] (G_3) at (0,0.5){};
  \node[vertex] (G_4) at (-0.5,1){};
  \node[vertex] (G_5) at (-0.5,1.5){};
  \node[vertex] (G_6) at (0.5,1.5){};
  \node[vertex] (G_7) at (0.5,1){};
  \node[vertex] (G_8) at (0,-0.5){};
  \node[vertex] (G_9) at (-0.5,-1){};
  \node[vertex] (G_10) at (-0.5,-1.5){};
  \node[vertex] (G_11) at (0.5,-1){};
  \node[vertex] (G_12) at (0.5,-1.5){};
  \draw[thick] (G_5)--(G_4) -- (G_3)--(G_8)--(G_9)--(G_10);
  \draw[thick] (G_6)--(G_7)--(G_3);
  \draw[thick] (G_12)--(G_11)--(G_8);
  \draw[thick] (G_1)--(G_3);
  \draw[thick] (G_2)--(G_8);
  \draw (0,1.2)node{$\ldots$};
  \draw (0,-1.2)node{$\ldots$};
  \draw (0,1.5)node{$r_7$};
  \draw (0,-1.5)node{$p_7$};
  \draw (0,-2)node{$T^7_{r_7,p_7}$};
  \end{tikzpicture}
  \hspace{3em}
  \begin{tikzpicture}[scale = 1.2]
  \tikzstyle{vertex}=[circle,fill=black,minimum size=0.38em,inner sep=0pt]
  \node[vertex] (G_1) at (-0.5,0.5){};
  \node[vertex] (G_2) at (-0.5,-0.5){};
  \node[vertex] (G_3) at (0,0.5){};
  \node[vertex] (G_4) at (-0.5,1){};
  \node[vertex] (G_5) at (-0.5,1.5){};
  \node[vertex] (G_6) at (0.5,1.5){};
  \node[vertex] (G_7) at (0.5,1){};
  \node[vertex] (G_8) at (0,-0.5){};
  \node[vertex] (G_9) at (-0.5,-1){};
  \node[vertex] (G_10) at (-0.5,-1.5){};
  \node[vertex] (G_11) at (0.5,-1){};
  \node[vertex] (G_12) at (0.5,-1.5){};
  \node[vertex] (G_13) at (0,0.17){};
  \node[vertex] (G_14) at (0,-0.17){};
  \draw[thick] (G_5)--(G_4) -- (G_3)--(G_13)--(G_14)--(G_8)--(G_9)--(G_10);
  \draw[thick] (G_6)--(G_7)--(G_3);
  \draw[thick] (G_12)--(G_11)--(G_8);
  \draw[thick] (G_1)--(G_3);
  \draw[thick] (G_2)--(G_8);
  \draw (0,1.2)node{$\ldots$};
  \draw (0,-1.2)node{$\ldots$};
  \draw (0,1.5)node{$r_8$};
  \draw (0,-1.5)node{$p_8$};
  \draw (0,-2)node{$T^8_{r_8,p_8}$};
  \end{tikzpicture}
    \caption{Trees $T^i_{r_i,p_i}$ for $3\leqslant i\leqslant 8.$}\label{fig04}
\end{figure}

Now we are ready to give the proof for Theorem~\ref{thm1.5}.
\begin{proof}[\bf Proof of Theorem~\ref{thm1.5}]
{\rm (i)} Let $G^*$ be a graph in $\mathcal{G}^2_n$ having the minimum $A_\alpha$-index. Then we know that $G^*$ does not contain $K_1\cup K_2$ and $3K_1$ as its induced subgraph. Let $\overline{G^*}$ be the complement of $G^*.$ Then we have $d_{\overline{G^*}}(v)\leqslant1$ for each vertex $v\in V_{G^*}$ (since otherwise there exist two vertices, say $u,w,$ of $V_{G^*}$ satisfying $uv,wv\notin E_{G^*},$ and thus $G^*[u,v,w]\cong K_1\cup K_2$ or $G^*[u,v,w]\cong 3K_1$). This implies that $E_{\overline{G^*}}$ is a matching of $K_n.$ Combining this with Lemma~\ref{lem2.1} gives us that $G^*\cong K_n-M$, where $M$ is a maximum matching of $K_n.$
\vspace{2mm}

{\rm (ii)} Since $\tau(P_n)=\left\lceil\frac{2}{3}n\right\rceil,$ the result follows immediately from Lemma~\ref{lem2.3}.
\vspace{2mm}

{\rm (iii)} Let $G^*$ be a graph in $\mathcal{G}^{n-1}_n$ having the minimum $A_\alpha$-index with $n\geqslant4.$ Since $n-1\geqslant\left\lceil\frac{2}{3}n\right\rceil,$ we have $G^*\in\mathscr{T}^{n-1}_n$ (based on Theorem~\ref{thm1.4}). Let $S=S'\cup S''$ be a maximum dissociation set of $G^*$ such that $S'$ is an independent set and $E_{G^*[S'']}$ is a matching of $G^*.$ Assume that $V_{G^*}=\{u_1,u_2,\ldots,u_n\}$ and $S=V_{G^*}\setminus\{u_1\}.$ Then $G^*\cong S_{k_1,k_2},$ where $S_{k_1,k_2}$ is depicted in Figure~\ref{fig01}, and $|S'|=k_1.$
\begin{claim}\label{C4}
$k_1\leqslant 1.$
\end{claim}
\begin{proof}[\bf Proof of Claim~\ref{C4}]
Suppose that $k_1\geqslant2$, i.e., there exist two vertices $u_s, u_t$ in $S',$ where $s,t\in\{2,\ldots,n\}.$ Then let $G^*_1=G^*-u_1u_s+u_tu_s.$ In view of Lemma~\ref{lem2.8}, we know that $G^*_1\in\mathscr{T}^{n-1}_n.$ However, by Lemma~\ref{lem2.5}, we have $\lambda_\alpha(G^*_1)<\lambda_\alpha(G^*),$ which contradicts the choice of $G^*.$ Therefore, $k_1\leqslant1,$ as desired.
\end{proof}
According to Claim~\ref{C4}, we obtain $G^*\cong S_{0,\frac{n-1}{2}}$ if $n$ is odd and $G^*\cong S_{1,\frac{n-2}{2}}$ if $n$ is even.
\vspace{2mm}

{\rm (iv)} Let $G^*$ be a graph in $\mathcal{G}_n^{n-2}$ having the minimum $A_\alpha$-index with $n\geqslant6.$
Since $n-2\geqslant\left\lceil\frac{2}{3}n\right\rceil,$ we have $G^*\in\mathscr{T}^{n-2}_n$ by Theorem~\ref{thm1.4}. Let $S=S'\cup S''$ be a maximum dissociation set of $G^*$ such that $S'$ is an independent set and $E_{G^*[S'']}$ is a matching of $G^*.$ Assume that $V_{G^*}=\{u_1,u_2,\ldots,u_n\}$ and $S=V_{G^*}\setminus\{u_1,u_2\},$ then we have $\min\{|N(u_1)\setminus\{u_2\}|,|N(u_2)\setminus\{u_1\}|\}\geqslant1.$

Note that $G^*$ is a tree. Then we know that $|N(u_1)\cap N(u_2)|\leqslant1,$ and thus there exist at least $|S'|-1$ leaves in $S'$. By a similar discussion as that in {\rm (iii)}, we obtain $|S'|\leqslant3.$ We first assume that $n$ is even. Since $n=|S|+2=|S'|+|S''|+2$ and $|S''|$ is even, we have $|S'|$ is even, i.e., $|S'|\in\{0,2\}.$ We split the proof in the following two cases.

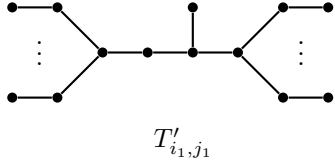
\begin{figure}[!ht]
\centering
  \begin{tikzpicture}[scale = 1.2]
  \tikzstyle{vertex}=[circle,fill=black,minimum size=0.38em,inner sep=0pt]
  \node[vertex] (G_1) at (0,0){};
  \node[vertex] (G_2) at (-0.5,0.5){};
  \node[vertex] (G_3) at (-1,0.5){};
  \node[vertex] (G_4) at (-0.5,-0.5){};
  \node[vertex] (G_5) at (-1,-0.5){};
  \node[vertex] (G_6) at (0.5,0){};
  \node[vertex] (G_7) at (1,0){};
  \node[vertex] (G_8) at (1,0.5){};
  \node[vertex] (G_9) at (1.5,0){};
  \node[vertex] (G_10) at (2,0.5){};
  \node[vertex] (G_11) at (2.5,0.5){};
  \node[vertex] (G_12) at (2,-0.5){};
  \node[vertex] (G_13) at (2.5,-0.5){};
  \draw[thick] (G_5)--(G_4) -- (G_1)--(G_6)--(G_7)--(G_9)--(G_10)--(G_11);
  \draw[thick] (G_3) -- (G_2)--(G_1);
  \draw[thick] (G_13) -- (G_12)--(G_9);
  \draw[thick] (G_8) -- (G_7);
  \draw (-0.7,0.08)node{$\vdots$};
  \draw (2.2,0.08)node{$\vdots$};
  \draw (0.9,-1)node{$T'_{i_1,j_1}$};
   \end{tikzpicture}
 \caption{Tree $T'_{i_1,j_1}$} \label{fig05}
\end{figure}

{\bf Case 1}\ \ $|N(u_1)\cap N(u_2)|=1.$ In this case, we proceed by considering the size of $S'$.

{\bf Subcase 1.1}\ \ $|S'|=0$. In this subcase, we have $G^*\cong T^3_{r_3,p_3}$ with $(r_3,p_3)\neq(0,0)$ and $r_3+p_3=\frac{n-4}{2}.$ If $\max\{r_3,p_3\}\leqslant1,$ then by Lemma~\ref{lem2.5}, we have $\lambda_\alpha(T^1_{r_3,p_3})<\lambda_\alpha(T^3_{r_3,p_3}).$ If $\max\{r_3,p_3\}\geqslant2,$ without loss of generality, we assume that $r_3\geqslant2.$
Let $T'_{i_1,j_1}$ be the tree obtained from $T^1_{i_1,j_1}$ by attaching a pendant edge to a 2-degree vertex on its internal path of length 3 (see Figure~\ref{fig05}). 
Then by Lemma~\ref{lem2.6}, we obtain
\begin{align}\label{eq:4.5}
\lambda_\alpha(T'_{r_3,p_3})<\lambda_\alpha(T^3_{r_3,p_3}).
\end{align}
It is easy to see that $T^1_{r_3,p_3}$ is a proper subgraph of $T'_{r_3,p_3},$ thus we have $\lambda_\alpha(T^1_{r_3,p_3})<\lambda_\alpha(T'_{r_3,p_3})$ by Lemma~\ref{lem2.1}. Together this with \eqref{eq:4.5}, we have $\lambda_\alpha(T^1_{r_3,p_3})<\lambda_\alpha(T^3_{r_3,p_3}).$ Note that $T^1_{r_3,p_3}\in\mathscr{T}^{n-2}_n.$ Both cases yield a contradiction to the choice of $G^*.$

{\bf Subcase 1.2}\ \ $|S'|=2$ and $N(u_1)\cap N(u_2)\subseteq S'$. In this subcase, we have $G^*\cong T^5_{r_5,p_5}$ with $r_5\geqslant0, p_5\geqslant1$ and $r_5+p_5=\frac{n-4}{2}.$ If $r_5=0,$ then $G^*\cong T^1_{0,p_5}.$ If $r_5\geqslant1, p_5=1,$ then by Lemma~\ref{lem2.5}, we have $\lambda_\alpha(T^1_{r_5,1})<\lambda_\alpha(G^*),$ a contradiction. If $r_5\geqslant1, p_5\geqslant2,$ then by Lemma~\ref{lem2.6}, we obtain
 \begin{align}\label{eq:4.6}
\lambda_\alpha(T^2_{r_5,p_5})<\lambda_\alpha(T^5_{r_5,p_5}).
\end{align}
Note that $T^1_{r_5,p_5}$ is a proper subgraph of $T^2_{r_5,p_5}.$ Therefore, by Lemma~\ref{lem2.1}, $\lambda_\alpha(T^1_{r_5,p_5})<\lambda_\alpha(T^2_{r_5,p_5}).$ Combining this with \eqref{eq:4.6} yields $\lambda_\alpha(T^1_{r_5,p_5})<\lambda_\alpha(T^5_{r_5,p_5}).$ Thus we derive a contradiction to the choice of $G^*$.

{\bf Subcase 1.3}\ \ $|S'|=2$ and $N(u_1)\cap N(u_2)\subseteq S''$. In this subcase, we have $G^*\cong T^6_{r_6,p_6}$ with $(r_6,p_6)\neq(0,0)$ and $r_6+p_6=\frac{n-6}{2}.$ Similarly, by Lemmas~\ref{lem2.1} and \ref{lem2.6}, we obtain $\lambda_\alpha(T^8_{r_6,p_6})<\lambda_\alpha(T^6_{r_6,p_6}).$ Note that $T^8_{r_6,p_6}\in\mathscr{T}^{n-1}_n,$ a contradiction again.

{\bf Case 2}\ \ $|N(u_1)\cap N(u_2)|=0.$ In this case, we proceed by considering the following four subcases.

{\bf Subcase 2.1}\ \ $|S'|=0$ and $u_1\sim u_2.$ In this subcase, we have $G^*\cong T^4_{r_4,p_4}$ with $\min\{r_4,p_4\}\geqslant1$ and $r_4+p_4=\frac{n-2}{2}.$ If $r_4=p_4=1,$ then $G^*\cong T^1_{1,0}.$ If $r_4=1$ and $p_4\geqslant2,$ then $G^*\cong T^1_{0,p_4}.$ If $r_4\geqslant2$ and $p_4=1$, then $G^*\cong T^1_{r_4,0}.$ If $r_4,p_4\geqslant2,$ then by Lemmas~\ref{lem2.1} and \ref{lem2.6}, we obtain
\begin{align}\label{eq:4.7}
\min\{\lambda_\alpha(T^1_{r_4-1,p_4}), \lambda_\alpha(T^1_{r_4,p_4-1})\}<\lambda_\alpha(T^1_{r_4,p_4})<\lambda_\alpha(G^*).
\end{align}
It is easy to see that both $T^1_{r_4-1,p_4}$ and $T^1_{r_4,p_4-1}$ are in $\mathscr{T}^{n-2}_n.$ Combining this with \eqref{eq:4.7} gives us a contradiction to the choice of $G^*.$

{\bf Subcase 2.2}\ \ $|S'|=0$ and $u_1\nsim u_2.$ In this subcase, we have $G^*\cong T^1_{r_1,p_1}$ with $r_1+p_1=\frac{n-4}{2}.$

{\bf Subcase 2.3}\ \ $|S'|=2$ and $u_1\sim u_2.$ In this subcase, we have $G^*\cong T^7_{r_7,p_7}$ with $\min\{r_7,p_7\}\geqslant1$ and $r_7+p_7=\frac{n-4}{2}.$ Similarly, by Lemmas~\ref{lem2.1} and \ref{lem2.6}, we get $\lambda_\alpha(T^1_{r_7,p_7})<\lambda_\alpha(T^8_{r_7,p_7})<\lambda_\alpha(G^*),$ and thus we get a contradiction.

{\bf Subcase 2.4}\ \ $|S'|=2$ and $u_1\nsim u_2.$ In this subcase, we have $G^*\cong T^8_{r_8,p_8}$ with $r_8+p_8=\frac{n-6}{2}.$ If $(r_8,p_8)=(0,0),$ then $G^*\cong T^1_{1,0},$ and the result holds. So, In what follows, we assume that $(r_8,p_8)\neq(0,0).$ If $r_8\geqslant1$ and $p_8=0,$ then by Lemma~\ref{lem2.5}, we obtain that $\lambda_\alpha(T^1_{r_8,1})<\lambda_\alpha(G^*).$ If $r_8=0$ and $p_8\geqslant1,$ then by Lemma~\ref{lem2.5} again, $\lambda_\alpha(T^1_{1,p_8})<\lambda_\alpha(G^*).$ Both cases give us a contradiction. Therefore, in what follows we may assume that $r_8,p_8\geqslant1$.
By a similar proof of Lemma~\ref{lem4.4}, applying Matlab to compare the spectral radii of the corresponding quotient matrices, we get $\min\{\lambda_\alpha(T^1_{r_8,p_8+1}),\lambda_\alpha(T^1_{r_8+1,p_8})\}<\lambda_\alpha(G^*),$ which leads to a contradiction since $\{T^1_{r_8,p_8+1}, T^1_{r_8+1,p_8}\}\subseteq \mathscr{T}^{n-2}_n.$

By Cases 1 and 2, we obtain that $G^*\cong T^1_{r_1,p_1}$ with $r_1+p_1=\frac{n-4}{2}.$ Without loss of generality, we may assume that $r_1\geqslant p_1.$ 
If $(r_1,p_1)=(1,0),$ then $G^*\cong T^1_{1,0}$, and thus the result holds. If $p_1=0$ and $r_1\geqslant2,$ then by Lemma~\ref{lem2.5}, we have $\lambda_\alpha(T^1_{r_1-1,1})<\lambda_\alpha(T^1_{r_1,0}),$ a contradiction. If $p_1\geqslant1,$ then in view of Lemma~\ref{lem4.4}, we obtain $G^*\cong T^1_{\left\lceil\frac{n-4}{4}\right\rceil,\left\lfloor\frac{n-4}{4}\right\rfloor},$ the result holds.
If $n$ is odd, we can also get $G^*\cong T^2_{\left\lfloor\frac{n-5}{4}\right\rfloor,\left\lceil\frac{n-5}{4}\right\rceil}$ in a similar way, whose procedure is omitted here.
\end{proof}
\section{\normalsize Concluding remarks}\setcounter{equation}{0}
In this paper, we mainly head in characterizing the extremal graphs having the maximum $A_\alpha$-index over all graphs in $\mathcal{G}_n^\tau,$  $\mathcal{B}_n^\tau$ and $\mathscr{T}_n^\tau$, respectively, and determining the structure of graphs having the minimum $A_\alpha$-index among all graphs in $\mathcal{G}_n^\tau$ with some restrictions on $\tau.$ Theorem~\ref{thm1.1} (resp. Theorem~\ref{thm1.2}, Theorem~\ref{thm1.3}) characterizes all the connected graphs (resp. bipartite graphs, trees) having the maximum $A_\alpha$-index among all connected graphs (resp. bipartite graphs, trees) with given order and dissociation number. Theorem~\ref{thm1.4} shows
that the graph over $\mathcal{G}_n^\tau$ having the minimum $A_\alpha$-index is a tree provided that $\tau\geqslant \left\lceil\frac{2}{3}n\right\rceil;$
whereas Theorem~\ref{thm1.5} determines the graphs with fixed order $n$ and dissociation number $\tau\in\{2,\left\lceil\frac{2}{3}n\right\rceil,n-1,n-2\}$ having the minimum $A_\alpha$-index.

In view of \eqref{eq:1.1}, if we put $\alpha=0$, respectively, in Theorems~\ref{thm1.1}--\ref{thm1.5}, then we may deduce the main results obtained for the adjacency spectral radius of graphs with given order and dissociation number (see \cite{HLZ} for details), whereas if we put $\alpha=\frac{1}{2}$, respectively, in Theorems~\ref{thm1.1}--\ref{thm1.5}, we may also deduce the corresponding results for the signless Laplacian spectral radius, say $q(G),$ of graphs with given order and dissociation number.  Then the next corollaries follow immediately.
\begin{cor}\label{cor5.1}
Let $G$ be in $\mathcal{G}_n^\tau$ having the maximum signless Laplacian spectral radius. Then $G\cong K_{n-\tau}\vee \left(\frac{\tau}{2}K_2\right)$ if $\tau$ is even, and $G\cong K_{n-\tau}\vee \left(\frac{\tau-1}{2}K_2\cup K_1\right)$ if $\tau$ is odd.
\end{cor}

\begin{cor}\label{cor5.2}
Let $G$ be a graph in $\mathcal{B}_n^\tau.$  Then $q(G)\leqslant n$ with equality if and only if $G\cong K_{\tau,n-\tau}$.
\end{cor}
\begin{cor}\label{cor5.3}
Let $T$ be a tree in $\mathscr{T}_n^\tau\ (n\geqslant 3)$ having the maximum signless Laplacian spectral radius. Then $q(T)\leqslant \theta(n,\tau)$
with equality if and only if $T\cong S^\dag_{n,\tau}$, where $\theta(n,\tau)$ is the largest root of
$x^3+(n-2\tau-6)x^2+(8\tau-4n+9)x-n=0.$
\end{cor}

\begin{cor}\label{cor5.4}
Let $G^\dag$ be a graph in $\mathcal{G}_n^\tau$ with $\tau\geqslant \left\lceil\frac{2}{3}n\right\rceil$ having the minimum signless Laplacian spectral radius. Then $G^\dag$ is a tree.
\end{cor}

\begin{cor}\label{cor5.5}
Let $G^*$ be a graph in $\mathcal{G}_n^\tau$ having the minimum signless Laplacian spectral radius.
\begin{wst}
\item[{\rm (i)}] If $\tau=2,$ then $G^*\cong K_n-M$, where $M$ is a maximum matching of $K_n.$
\item[{\rm (ii)}] If $\tau=\left\lceil\frac{2}{3}n\right\rceil,$ then $G^*\cong P_n.$
\item[{\rm (iii)}] If $n\geqslant4$ and $\tau=n-1,$ then $G^*\cong S_{0,\frac{n-1}{2}}$ if $n$ is odd and $G^*\cong S_{1,\frac{n-2}{2}}$ if $n$ is even.
\item[{\rm (iv)}] If $n\geqslant6$  and $\tau=n-2,$ then $G^*\cong T^1_{\left\lceil\frac{n-4}{4}\right\rceil,\left\lfloor\frac{n-4}{4}\right\rfloor}$ if $n$ is even and $G^*\cong T^2_{\left\lfloor\frac{n-5}{4}\right\rfloor,\left\lceil\frac{n-5}{4}\right\rceil}$ if $n$ is odd.
\end{wst}
\end{cor}



\section*{\normalsize Declaration of competing interest}
There is no competing interest.

\appendix
\renewcommand{\thesection}{\normalsize Appendix}
\section{\normalsize}
\renewcommand{\thesection}{\normalsize \Alph{section}}
\begin{lstlisting}[language=Matlab]
clear
clc
syms a b c alpha x
AA=[alpha*(b+1) 0 (1-alpha)*b 1-alpha
     0 alpha*(a+b) (1-alpha)*b (1-alpha)*a
     (1-alpha)*a (1-alpha)*c alpha*(a+c) 0
     1-alpha (1-alpha)*c 0 alpha*(c+1)];
I=eye(4);
y=simplify(det(x*I-AA))
fun=matlabFunction(y);
objfun=@(x)fun(x(1),x(2),x(3),x(4),x(5));
x0=[0 0 0 0 0];
A=[-1 0 0 1 0;0 0 1 -1 0];
b=[-1 0];
Aeq=[];
beq=[];
lb=[2 0 0 1 0];
ub=[inf 1 inf inf inf];
[x,fval]=fmincon(objfun,x0,A,b,Aeq,beq,lb,ub,@confun)
function [c,ceq]=confun(x)
c=0.5.*(x(2).*(2.*x(1)+x(3)+x(4))+sqrt(x(2).^2.*(2.*x(1)+x(3)+x(4)).^2
+8.*x(1).*(x(3)+x(4)).*(1-2.*x(2))))-x(5);
ceq=0;
end
\end{lstlisting}


\begin{thebibliography}{99}
\small \setlength{\itemsep}{-.8mm}
\bibitem{ABKL} V.E. Alekseev, R. Boliac, D.V. Korobitsyn, V.V. Lozin, NP-hard graph problems and boundary classes of graphs, Theoret. Comput. Sci. 389 (1-2) (2007) 219-236.
\bibitem{BPPR} F. Bock, J. Pardey, L. Penso, D. Rautenbach, A bound on the dissociation number, J. Graph Theory,  2023;103:661-673.
\bibitem{BPPR1} F. Bock, J. Pardey, L. Penso, D. Rautenbach, Relating the independence number and the dissociation number, J. Graph Theory. 1-21 (2023) DOI: 10.1002/jgt.22965.
\bibitem{BCL} R. Boliac, K. Cameron, V.V. Lozin, On computing the dissociation number and the induced matching number of bipartite graphs, Ars Combin. 72 (2004) 241-253.
\bibitem{CH} K. Cameron, P. Hell, Independent packings in structured graphs, Math. Program. 105 (2006) 201-213.
\bibitem{CLM} Y.Y. Chen, D. Li, and J.X. Meng, On the second largest $A_\alpha$-eigenvalues of graphs, Linear Algebra Appl. 580 (2019) 343-358.
\bibitem{DM2023}J. Das, S. Mohanty, Maximization of the spectral radius of block graphs with a given dissociation number, arXiv:2301.12790, math.CO
\bibitem{D-R} M. Desai, V. Rao, A characterization of the smallest eigenvalue of a graph, J. Graph Theory 18 (2) (1994) 181-194.
\bibitem{G-R} C. Godsil, G. Royle, Algebraic Graph Theory, vol. 207 of Graduate Texts in Mathematics, Springer-Verlag, New York, 2001.
\bibitem{CS} H.Y. Guo, B. Zhou, On the $\alpha$-spectral radius of graphs, Appl. Anal. Discrete Math. 14 (2) (2020) 431-458.
\bibitem{HLZ} J. Huang, X.Y. Geng, S.C. Li, Z.H. Zhou, On spectral extrema of graphs with given order and dissociation number, Discrete Appl. Math. 342 (2024) 368-380.
\bibitem{HLX} X. Huang, H.Q. Lin, J. Xue, The Nordhaus-Gaddum type inequalities of $A_\alpha$-matrix, Appl. Math. Comput. 365 (2020) 124716.
\bibitem{FLZ} D. Li, Y.Y. Chen, J.X. Meng, The $A_{\alpha}$-spectral radius of trees and unicyclic graphs with given degree sequence, Appl. Math. Comput. 363 (2019) 124622.
\bibitem{LS} W.T. Sun, S.C. Li, On the maximal number of maximum dissociation sets in forests with fixed order and dissociation number, Taiwanese J. Math. 1-37
(2023) DOI: 10.11650/tjm/230204.
\bibitem{L-S} S.C. Li, W.T. Sun, An arithmetic criterion for graphs being determined by their generalized $A_\alpha$-spectra, Discrete Math. 344 (8) (2021) 112469.
\bibitem{LTYZ} S.C. Li, B-S. Tam, Y.T. Yu, Q. Zhao, Connected graphs of fixed order and size with maximal $A_\alpha$-index: the one-dominating-vertex case, Linear Algebra Appl. 662 (2023) 110-135.
\bibitem{L-W} S.C. Li, S.J. Wang, The $A_\alpha$-spectrum of graph product, Electron. J. Linear Algebra 35 (2019) 473-481.
\bibitem{LY-2}S.C. Li, Y.T. Yu, On $A_\alpha$ spectral extrema of graphs forbidding even cycles, Linear Algebra Appl. 668 (2023), 11-27.
\bibitem{LY-1}S.C. Li, Y.T. Yu, H.H. Zhang, An $A_\alpha$-spectral Erd\H{s}-P\'osa theorem, Discrete Math. 346 (9) (2023) 113494.
\bibitem{L-Z} S.C. Li, Z.H. Zhou, On the $A_\sigma$-spectral radii of graphs with some given parameters, Rocky Mountain J. Math. 52 (3) (2022) 949-966.
\bibitem{LHX} H.Q. Lin, X. Huang, and J. Xue, A note on the $A_\alpha$-spectral radius of graphs, Linear Algebra Appl. 557 (2018) 430-437.
\bibitem{V.N} V. Nikiforov, Merging the $A$- and $Q$-spectral theories, Appl. Anal. Discrete Math. 11 (2017) 81-107.
\bibitem{C.Y} V. Nikiforov, G. Pastn, O. Rojo, R. L. Soto, On the $A_{\alpha}$-spectra of trees, Linear Algebra Appl. 520 (2017) 286-305.
\bibitem{N-O} V. Nikiforov and O. Rojo, On the $\alpha$-index of graphs with pendent paths, Linear Algebra Appl. 550 (2018) 87-104.
\bibitem{OD} Y. Orlovich, A. Dolgui, G. Finke, V. Gordon, F. Werner, The complexity of dissociation set problems in graphs, Discrete Appl. Math. 159 (13) (2011) 1352-1366.
\bibitem{PY} C. H. Papadimitriou, M. Yannakakis, The complexity of restricted spanning tree problems, J. Assoc. Comput. Mach. 29 (2) (1982) 285-309.
\bibitem{TLD}J.H. Tu, Y.X. Li, J.F. Du, Maximal and maximum dissociation sets in general and triangle-free graphs, Appl. Math. Comput. 426 (2022) 127107.
\bibitem{TZD} J.H. Tu, L. Zhang, J.F. Du, On the maximum number of maximum dissociation  sets in trees with given dissociation number, arXiv:2103.01407v1.
\bibitem{TZS} J.H. Tu, Z.P. Zhang, Y.T. Shi, The maximum number of maximum dissociation sets in trees, J. Graph Theory 96 (4) (2021) 472-489.
\bibitem{WWT} S. Wang, D. Wong, F.L. Tian, Bounds for the largest and the smallest $A_\alpha$ eigenvalues of a graph in terms of vertex degrees, Linear Algebra Appl. 590 (2020) 210-223.
\bibitem{West} D.B. West, Introduction to Graph Theory, Prentice Hall, Inc., Upper Saddle River, NJ, 1996.
\bibitem{XWT} F. Xu, D. Wong, F.L. Tian, On the multiplicity of $\alpha$ as an eigenvalue of the $A_\alpha$ matrix of a graph in terms of the number of pendant vertices, Linear Algebra Appl. 594 (2020) 193-204.
\bibitem{HZ} J. Xue, H.Q. Lin, S.T. Liu, and J.L. Shu, On the $A_{\alpha}$-spectral radius of a graph, Linear Algebra Appl. 550 (2018) 105-120.
\bibitem{Y81} M. Yannakakis, Node-deletion problems on bipartite graphs, SIAM J. Comput. 10 (2) (1981) 310-327.
\bibitem{YYSX} L.H. You, M. Yang, W. So, W.G. Xi, On the spectrum of an equitable quotient matrix and its application, Linear Algebra Appl. 577 (2019) 21-40.
\end{thebibliography}
\end{document}